\documentclass[11pt]{amsart}
\usepackage{amsmath, amssymb, amsthm, amsfonts}
\usepackage{enumerate}
\usepackage{dsfont}
\usepackage{cite}
\usepackage[colorlinks, citecolor=red]{hyperref}
\usepackage{mathrsfs}
\usepackage{epsfig}
\usepackage{lscape}
\usepackage{subfigure}
\usepackage{epstopdf}
\usepackage{caption}
\usepackage{algorithm}
\usepackage{algpseudocode}
\usepackage{multicol}  
\usepackage{multirow}
\usepackage{geometry}
\usepackage{paralist}
\usepackage{url}
\usepackage{graphicx}
\usepackage{hyperref}
\usepackage{booktabs}
\usepackage{tcolorbox}
\usepackage{setspace} 
\usepackage[numbers]{natbib}
\usepackage{longtable}
\usepackage{tikz}
\usepackage{enumitem}
\usepackage{array}

\usepackage{tikz-3dplot}
\usepackage{amsmath}
\usepackage{placeins}
\usetikzlibrary{3d, arrows.meta, calc, decorations.pathreplacing}
\usetikzlibrary{positioning,decorations.markings}

\theoremstyle{definition}
\newtheorem{definition}{Definition}[section]

\newtheorem{theorem}[definition]{Theorem}

\newtheorem{remark}[definition]{Remark}

\newtheorem{proposition}{Proposition}

\date{}

\tikzset{
	vector/.style={->, thick, >=Latex}, 
	subspace/.style={draw, dashed, rounded corners=5pt, opacity=0.3}, 
	ortho/.style={decorate, decoration={markings, mark=at position 0.5 with {\node[transform shape] {$\perp$};}}} 
}
 
\begin{document}
	\baselineskip 18pt
	\bibliographystyle{abbrv}
	\title{Projecting onto the unit dual quaternion set}
	
	\author{Ziyang Li}
	\address{School of Mathematical Sciences, Beihang University, Beijing, 100191, China. }
	\email{lzyang0114@buaa.edu.cn}
	
	\author{Chunfeng Cui}
	\address{LMIB of the Ministry of Education, School of Mathematical Sciences, Beihang University, Beijing, 100191, China. }
	\email{chunfengcui@buaa.edu.cn}
	
	\author{Jiaxin Xie}
	\address{LMIB of the Ministry of Education, School of Mathematical Sciences, Beihang University, Beijing, 100191, China. }
	\email{xiejx@buaa.edu.cn}
	\begin{abstract}
		Dual quaternions have gained significant attention due to their wide applications in areas such as multi-agent formation control, 3D motion modeling, and robotics. A fundamental aspect in dual quaternion research involves the projection onto the unit dual quaternion set. In this paper, we systematically study such projections under the $2^R$-norm, which is commonly used in practical applications. We identify several distinct cases based on the relationship between the standard and dual parts in vector form, and demonstrate the effectiveness of the proposed algorithm through numerical experiments.
	\end{abstract}
	
	\maketitle

	\let\thefootnote\relax
	
	\let\thefootnote\relax

	\section{Introduction}
	
	
	Dual quaternions, introduced by Clifford \cite{Clifford1871PreliminarySO} in 1871, form an 8-dimensional algebra combining quaternions and dual numbers that elegantly represents the space of 3D rigid transformations \cite{brambley2020unit,han2008kinematic,  wang2012dual,giribet2025dual}. This mathematical framework has become indispensable across engineering fields including robotics, computer graphics, and biomechanics due to its superior handling of combined rotational and translational motions \cite{cheng2016dual,daniilidis1996dual,daniilidis1999hand,wang2024dynamic}. This superiority is particularly evident in its ability to unify general screw motions with all fundamental cases of rotation and translation, such as the rotation of a robotic joint or the linear movement of a camera on a track.
	
	Among the many research problems involving dual quaternions, projection onto the unit dual quaternion set is of fundamental importance. For instance, the hand-eye calibration problem can be formulated as a least-squares optimization of a real-valued function over dual quaternion variables, subject to a unit dual quaternion constraint \cite{Can2023APL,qi2023standard,tao2024solving,wang2024qly,ling2025metric}. When solving such problems using a proximal linearized algorithm, each iteration requires computing a projection onto the unit dual quaternion set \cite{Regularization-Patching,zhu2025subspace,cui2024power}. 
	The set of unit dual quaternions forms a nonlinear constrained set with intrinsic geometric structure. Recent work models it as a Riemannian manifold and performs pose estimation and optimization directly on the manifold, where the unit constraint is preserved intrinsically rather than enforced via explicit projections \cite{li2020improved,warke2024pose}.

	In this paper, we aim to systematically investigate the projection onto the unit dual quaternion set under the \(2^R\)-norm, a norm widely adopted in practical applications as it ensures that the rotational and translational components are treated with equal importance. 
	In particular, 
	we develop a unified projection framework based on the first-order optimality condition (KKT conditions) for constrained optimization, with particular attention to the distinct characteristics between the standard and dual parts of dual quaternions. 
	Within this framework, we provide a complete case classification of the projection problem and analyze each case in a systematic manner.
	For several special cases, we further derive closed-form expressions of the projection points. Compared with the work in  \cite{Can2023APL}, our analysis provides a more comprehensive and systematic treatment of the projection problem; further discussion is provided in Remark \ref{remark-2}.
	
	The rest of the paper is organized as follows. Section \ref{SEC2} reviews the fundamentals of quaternions and dual quaternions, then introduces the first-order optimality conditions for constrained optimization. Section \ref{SEC3} proposes the algorithm for projecting onto the unit dual quaternion set, which comprehensively accounts for all four possible cases. Section \ref{SEC4} presents numerical examples of different standard-dual part configurations. Finally, some conclusions are drawn in Section \ref{SEC5}.
	
	\section{Preliminaries}\label{SEC2}
	In this section, we review the definitions of quaternions, dual numbers, and dual quaternions. The sets of real numbers, quaternions, dual numbers, dual quaternions, and unit dual quaternions are denoted as $\mathbb{R}, \mathbb{Q}, \hat{\mathbb{R}}, \hat{\mathbb{Q}}$, and $\hat{\mathbb{U}}$ respectively. 
	
	\subsection{Quaternions}
	
	A quaternion \(\tilde{q} \in \mathbb{Q}\) is defined as a pair consisting of a scalar part \(q_0 \in \mathbb{R}\) and a vector part \(\overrightarrow{q} = (q_1, q_2, q_3)^\top \in \mathbb{R}^3\). To distinguish quaternions from real numbers, we adopt the convention of using a tilde symbol. A quaternion can then be written in either of the following equivalent forms:
	\begin{equation*}
		\tilde{q} = [q_0, \overrightarrow{q}] = [q_0, q_1, q_2, q_3] = q_0 + \mathbf{i}q_1 + \mathbf{j}q_2 + \mathbf{k}q_3,
	\end{equation*}
	where \(q_0, q_1, q_2, q_3 \in \mathbb{R}\), and \(\mathbf{i}, \mathbf{j}, \mathbf{k}\) are the standard quaternion imaginary units, which satisfy
	\[
	\mathbf{i}^2 = \mathbf{j}^2 = \mathbf{k}^2 = \mathbf{ijk} = -1,\quad \mathbf{ij=-ji=k,\ ki=-ik=j}, \ \text{and} \
	\mathbf{jk=-kj=i}.
	\]
	If $q_0=0$, then $\tilde{q}=[0, \overrightarrow{q}]$ is called a pure quaternion.
	The conjugate of \(\tilde{q}\), denoted \(\tilde{q}^*\), is defined by
	$
	\tilde{q}^* = q_0 - \mathbf{i}q_1 - \mathbf{j}q_2 - \mathbf{k}q_3.
	$
	The scalar part of \(\tilde{q}\), denoted by \(\operatorname{Sc}(\tilde{q})\), is given by
	$
	\operatorname{Sc}(\tilde{q}) = \frac{1}{2} (\tilde{q} + \tilde{q}^*) = q_0.
	$
	Given two quaternions \(\tilde{p} = [p_0, \vec{p}]\) and \(\tilde{q} = [q_0, \vec{q}]\), their sum and product are defined as follows:
	\[	
	\tilde{q}\pm \tilde{p}=\left[q_0\pm p_0,\overrightarrow{q}\pm\overrightarrow{p}\right] \  \text{and} \ \tilde{q} \tilde{p}=\left[ q_0p_0-\overrightarrow{q}\cdot \overrightarrow{p},\,p_0\overrightarrow{q}+q_0\overrightarrow{p}+\overrightarrow{q}\times\overrightarrow{p}\right] .
	\]
	Quaternion multiplication is associative but not commutative in general. That is, \(\tilde{q} \tilde{p} \neq \tilde{p} \tilde{q}\), unless the cross product \(\vec{q} \times \vec{p} = \vec{0}\). The magnitude (or norm) of \(\tilde{q}\) is defined as
	$
	|\tilde{q}| = \sqrt{\tilde{q}\tilde{q}^*} = \sqrt{\tilde{q}^* \tilde{q}} = \sqrt{q_0^2 + q_1^2 + q_2^2 + q_3^2}.
	$
	A quaternion \(\tilde{q}\) is called a unit quaternion if \(|\tilde{q}|^2 = 1\).
	
	\subsection{Dual quaternions}\label{Dual quaternions}
	Before introducing dual quaternions, we first present the concept of dual numbers. A dual number \( a = a_s + a_d \epsilon \in \hat{\mathbb{R}} \) consists of a standard part \( a_s \in \mathbb{R} \) and a dual part \( a_d \in \mathbb{R} \). The symbol \( \epsilon \) is called the infinitesimal unit, satisfying \( \epsilon \neq 0 \), \( \epsilon^2 = 0 \). 
	If \( a_s \neq 0 \), then \( a \) is said to be appreciable; otherwise, it is called infinitesimal~\cite{qi2022dual}. For two dual numbers \( a = a_s + a_d \epsilon \) and \( b = b_s + b_d \epsilon \) in \( \hat{\mathbb{R}} \), their addition, multiplication, and division are defined as:
	\begin{align*}
		a + b &= a_s + b_s + (a_d + b_d)\epsilon, \\
		ab &= a_s b_s + (a_s b_d + b_s a_d)\epsilon, \\
		\frac{b_s + b_d \epsilon}{a_s + a_d \epsilon} &=
		\begin{cases}
			\displaystyle\frac{b_s}{a_s} + \left( \frac{b_d}{a_s} - \frac{b_s a_d}{a_s^2} \right)\epsilon, & \text{if } a_s \neq 0,\\[1ex]
			\displaystyle\frac{b_d}{a_d} + c\epsilon, & \text{if } a_s = 0,\, b_s = 0,
		\end{cases}
	\end{align*}
	where \( c \in \mathbb{R} \) is an arbitrary real constant. 
	
	Let the set of dual quaternions be defined as
	\begin{equation*}
		\hat{\mathbb{Q}} := \left\{ \hat{q} = \tilde{q}_s + \tilde{q}_d \epsilon \,\middle|\, \tilde{q}_s, \tilde{q}_d \in \mathbb{Q} \right\}.
	\end{equation*}
	We use a hat symbol ``\(\:\hat{}\:\)'' to distinguish dual quaternions from dual numbers.
	In a dual quaternion \(\hat{q} = \tilde{q}_s + \tilde{q}_d \epsilon\), the component \(\tilde{q}_s\) is called the standard part, and \(\tilde{q}_d\) is the dual part. Addition and subtraction of dual quaternions follow the same rules as for dual numbers. The conjugate of a dual quaternion is defined component-wise as
	$
	\hat{q}^* = \tilde{q}_s^* + \tilde{q}_d^* \epsilon.
	$
	The magnitude of \(\hat{q}\) ~\cite{qi2022dual} is defined as
	\begin{equation}\label{magnitude of DQ}
		|\hat{q}| =
		\begin{cases}
			|\tilde{q}_s| + \dfrac{\operatorname{Sc}(\tilde{q}_s^* \tilde{q}_d)}{|\tilde{q}_s|} \epsilon, & \text{if } \tilde{q}_s \neq 0, \\[1ex]
			|\tilde{q}_d| \, \epsilon, & \text{otherwise}.
		\end{cases}
	\end{equation}
	A dual quaternion \(\hat{q} = \tilde{q}_s + \tilde{q}_d \epsilon\) is called a unit dual quaternion if it satisfies the following two conditions:
	\begin{equation}\label{unit dual quaternion}
		|\tilde{q}_s| = 1 \quad \text{and} \quad \tilde{q}_s^* \tilde{q}_d + \tilde{q}_d^* \tilde{q}_s = 0.
	\end{equation}
	The magnitude of a dual quaternion is associated with its real and dual parts. As a result, the unit dual quaternion constraint consists of two coupled conditions: the real part must have unit norm, and the real and dual parts must satisfy an orthogonality condition.
	For a dual quaternion \(\hat{q} = \tilde{q}_s + \tilde{q}_d \epsilon \in \hat{\mathbb{Q}}\), the \(2^R\)-norm is defined as
	\begin{equation*}
		\| \hat{q} \|_{2^{R}} := \sqrt{|\tilde{q}_s|^2 + |\tilde{q}_d|^2}.
	\end{equation*}
	
	Unit dual quaternions offer a canonical and compact representation of rigid body motion. A single unit dual quaternion $\hat{q} = \tilde{q}_s + \tilde{q}_d\epsilon$, subject to the two conditions in \eqref{unit dual quaternion}, simultaneously encodes both rotation and translation. The rotation matrix $\mathbf{R}=(r_{ij})\in \mathbb{R}^{3\times 3}$ is represented by the unit quaternion $\tilde{q}_s=q_{s_0}+q_{s_1}\mathbf{i}+q_{s_2} \mathbf{j}+q_{s_3} \mathbf{k}$, where $q_{s_0}=\frac{1}{2}\sqrt{1+r_{11}+r_{22}+r_{33}}$, $q_{s_1}=\frac{1}{4q_{s_0}}(r_{23}-r_{32})$, $q_{s_2}=\frac{1}{4q_{s_0}}(r_{13}-r_{31})$ and $q_{s_3}=\frac{1}{4q_{s_0}}(r_{12}-r_{21})$. When we know the unit quaternion $\tilde{q}_s$, we have the rotation matrix:
	\begin{equation*}
		\mathbf{R}=
		\begin{pmatrix}
			q_{s_0}^2+q_{s_1}^2-q_{s_2}^2-q_{s_3}^2 & 2(q_{s_1}q_{s_2}-q_{s_0}q_{s_3}) & 2(q_{s_1}q_{s_3}+q_{s_0}q_{s_2})\\
			2(q_{s_1}q_{s_2}+q_{s_0}q_{s_3}) & q_{s_0}^2-q_{s_1}^2+q_{s_2}^2-q_{s_3}^2 & 2(q_{s_2}q_{s_3}-q_{s_0}q_{s_1})\\
			2(q_{s_1}q_{s_3}-q_{s_0}q_{s_2}) & 2(q_{s_2}q_{s_3}+q_{s_0}q_{s_1}) & q_{s_0}^2-q_{s_1}^2-q_{s_2}^2+q_{s_3}^2
		\end{pmatrix}.
	\end{equation*}
	The translation vector $\mathbf{t}=(t_1,t_2,t_3)^\top$ is encoded within the dual part $\tilde{q}_d$, and is recovered as a pure quaternion $\tilde{t}=0+t_1\mathbf{i}+t_2\mathbf{j}+t_3\mathbf{k}$ using the relation $\tilde{t}=2\tilde{q}_d\tilde{q}_s^*$. One may refer to \cite{vince2011rotation} for more information on the relationship between quaternions and rotations.
	
	\subsection{First-order optimality conditions}
	Consider the following constrained optimization problem
	\begin{equation}
		\label{const-opt}
		\min_{x \in \mathbb{R}^n } f(x) \quad \text{subject to} \quad
		\begin{cases}
			c_i(x) = 0, & i \in \mathcal{E}, \\
			c_i(x) \geq 0, & i \in \mathcal{I},
		\end{cases}
	\end{equation}
	where $f$ and the functions $c_i$ are all smooth, real-valued functions on a subset of $\mathbb{R}^n$, and $\mathcal{I}$ and $\mathcal{E}$ are two finite sets of indices. We call $f$ the objective function, while $c_i(x)=0$, $i\in \mathcal{E}$ are the equality constraints and $c_i(x)\geq0$, $i \in \mathcal{I}$ are the inequality constraints. 
	
	\begin{definition}[Active set] The active set $\mathcal{A}(x)$ at any feasible $x$ consists of the equality constraint indices from $\mathcal{E}$ together with the indices of the inequality constraints $i$ for which $c_i(\textbf{x})=0$, that is,
		\begin{equation*}
			\mathcal{A}(x)=\mathcal{E}\cup \left\lbrace  i\in \mathcal{I}|c_i(x)=0\right\rbrace.
		\end{equation*}
	\end{definition}
	
	\begin{definition}[LICQ] Given the point $x$ and the active set $\mathcal{A}(x)$, we say that the linear independence constraint qualification (LICQ) holds if the set of active constraint gradients $\left\lbrace \nabla c_i(x), i\in \mathcal{A}(x)\right\rbrace $ is linearly independent.
	\end{definition}
	
	Associated with each constraint $c_i(x)$ in \eqref{const-opt}, we introduce a Lagrange multiplier $\lambda_i \in \mathbb{R}$. Then we define the Lagrangian function for the general problem \eqref{const-opt} as
	$$
	\mathcal{L}(x,\mathbf{\lambda})=f(x)-\sum_{i\in\mathcal{E}\cup\mathcal{I}}\lambda_i c_i(x),
	$$
	where $\lambda$ is the vector of multipliers $\lambda_i$.
	
	\begin{theorem}[First-order necessary conditions] \label{First-Order}
		Suppose that $x^\# \in \mathbb{R}^n$ is a local solution of \eqref{const-opt}, the functions $f$ and $c_i$ in \eqref{const-opt} are continuously differentiable, and the LICQ holds at $x^\#$. Then there is a Lagrange multiplier vector $\lambda^\#\in \mathbb{R}^{|\mathcal{E}| +|\mathcal{I}|}$, with components $\lambda_i^\#$, $i \in \mathcal{E}\cup \mathcal{I}$, such that the following conditions are satisfied at $(x^\#,\lambda^\#)$
		\[	\begin{aligned}
			\nabla_{x}\mathcal{L}(x^\#,\lambda^\#)=0,&\\
			c_i(x^\#)=0,& \ \ \text{for all} \ i\in\mathcal{E},\\
			c_i(x^\#)\geq 0,& \ \ \text{for all} \ i \in \mathcal{I},\\
			\lambda_i^\#\geq0,& \ \ \text{for all} \ i \in \mathcal{I},\\
			\lambda_i^\#c_i(x^\#)=0,& \ \ \text{for all} \ i\in\mathcal{I}.
		\end{aligned}
		\]
	\end{theorem}
	The conditions above are widely known as the \textit{Karush-Kuhn-Tucker conditions}, or \textit{KKT conditions} for short \cite{nocedal2006numerical}.

	\section{Projection onto the unit dual quaternion set}\label{SEC3}
	\subsection{Problem setup}
	Suppose we have a dual quaternion \(\hat{a}\). Its projection onto the unit dual quaternion set under the \(2^R\)-norm can be formulated as the following optimization problem:
	\begin{equation}\label{problem-dq}
		\min_{\hat{q}} \:\frac{1}{2} \left\|  \hat{q} - \hat{a} \right\| _{2^R}^2 \quad \text{subject to} \quad |\hat{q}| = 1.
	\end{equation}
	To proceed, we introduce the real vector representation of dual quaternions and reformulate problem~\eqref{problem-dq} into a non-convex optimization problem. To distinguish between scalars and vectors, we adopt the convention that plain symbols such as \(0\) denote scalars, while arrows denote vectors. For example, \(\overrightarrow{0} = (0, 0, 0, 0)^\top\) denotes a four-dimensional zero vector.
	
	A dual quaternion \(\hat{q} = \tilde{q}_s + \tilde{q}_d \epsilon\) can be represented in real vector form as \(\overrightarrow{\hat{q}}\), namely
	\begin{equation*}
		\overrightarrow{\hat{q}} = \begin{pmatrix} \overrightarrow{\tilde{q}_s}\\ \overrightarrow{\tilde{q}_d} \end{pmatrix} =
		(q_{s_0}, q_{s_1}, q_{s_2}, q_{s_3}, q_{d_0}, q_{d_1}, q_{d_2}, q_{d_3})^\top, \quad
		\overrightarrow{\hat{q}} \in \mathbb{R}^8,\quad \overrightarrow{\tilde{q}_s},\, \overrightarrow{\tilde{q}_d } \in \mathbb{R}^4,
	\end{equation*}
	where \(\overrightarrow{\tilde{q}_s}= (q_{s_0}, q_{s_1}, q_{s_2}, q_{s_3})^\top\) and \(\overrightarrow{\tilde{q}_d} = (q_{d_0}, q_{d_1}, q_{d_2}, q_{d_3})^\top\) denote the vectorized real and dual parts of \(\hat{q}\), respectively. Thus, \(\overrightarrow{\hat{q}}\) can be viewed as an element in an eight-dimensional real vector space spanned by the standard basis \((1, \mathbf{i}, \mathbf{j}, \mathbf{k}, \epsilon, \epsilon \mathbf{i}, \epsilon \mathbf{j}, \epsilon \mathbf{k})\).
	Using this representation, we now translate the unit dual quaternion constraint \(|\hat{q}| = 1\) into constraints on \(\overrightarrow{\tilde{q}_s}\) and \(\overrightarrow{\tilde{q}_d }\). From the unit dual quaternion condition \eqref{unit dual quaternion}, it is known that
	\[
	\left\| \overrightarrow{\tilde{q}_s}\right\| _2^2 = 1
	\]
	and
	\begin{align*}
		0 = &\tilde{q}_s^* \tilde{q}_d + \tilde{q}_d^* \tilde{q}_s \\
		= &\left[q_{s_0} q_{d_0} + \overrightarrow{q_s} \cdot \overrightarrow{q_d},\; -q_{d_0} \overrightarrow{q_s} + q_{s_0} \overrightarrow{q_d} - \overrightarrow{q_s} \times \overrightarrow{q_d}\right] \\
		& + \left[q_{s_0} q_{d_0} + \overrightarrow{q_s} \cdot \overrightarrow{q_d},\; q_{d_0} \overrightarrow{q_s} - q_{s_0} \overrightarrow{q_d} - \overrightarrow{q_d} \times \overrightarrow{q_s}\right] \\
		= &\left[2\left(q_{s_0} q_{d_0} + \overrightarrow{q_s} \cdot \overrightarrow{q_d}\right), \; \overrightarrow{0} \right]\\
		=&\left[2\overrightarrow{\tilde{q}_d }^\top \overrightarrow{\tilde{q}_s },\ \overrightarrow{0}\right],
	\end{align*}
	which implies that \(\overrightarrow{\tilde{q}_d }^\top \overrightarrow{\tilde{q}_s }= 0\).
	Hence, in real vector form, the unit dual quaternion set corresponds to the following set
	\begin{equation*}
		\Omega := \left\{ \left( \overrightarrow{\tilde{q}_s}, \overrightarrow{\tilde{q}_d }\right)  \in \mathbb{R}^8 \ \middle| \
		\left\|  \overrightarrow{\tilde{q}_s }\right\| _2^2 = 1 \ \text{and} \ \overrightarrow{\tilde{q}_d}^\top \overrightarrow{\tilde{q}_s }= 0 \right\}.
	\end{equation*}
	Therefore, using the real vector representation, the original projection problem~\eqref{problem-dq} can be reformulated as the following non-convex optimization problem
	\begin{equation}\label{Objective Function}
		\begin{aligned}
			\min_{\overrightarrow{\tilde{q}_s},\, \overrightarrow{\tilde{q}_d}}& \quad f\left( \overrightarrow{\tilde{q}_s}, \overrightarrow{\tilde{q}_d}\right)  :=
			\frac{1}{2} \left\| \overrightarrow{\tilde{q}_s} - \overrightarrow{\tilde{a}_s}\right\| _2^2 + \frac{1}{2} \left\| \overrightarrow{\tilde{q}_d} - \overrightarrow{\tilde{a}_d}\right\| _2^2 \\
			\text{subject to}& \quad \left\| \overrightarrow{\tilde{q}_s}\right\| _2^2 = 1, \quad \overrightarrow{\tilde{q}_d}^\top \overrightarrow{\tilde{q}_s} = 0.
		\end{aligned}
	\end{equation}

	\subsection{The existence of the global minimum of \eqref{Objective Function}} 
	
	In this subsection, we show that problem~\eqref{Objective Function} admits a global minimizer.
	We first recall the extreme value theorem.
	\begin{theorem}[Extreme value theorem {\cite[Theorem~4.16]{RudinPMA}}]\label{EVT}
		Suppose $f$ is a continuous real-valued function on a compact metric space $\Lambda$, and
		\begin{equation}\label{eq:EVT-sup-inf}
			M=\sup_{p\in \Lambda} f(p), \qquad m=\inf_{p\in \Lambda} f(p).
		\end{equation}
		Then there exist points $p,q\in \Lambda$ such that $f(p)=M$ and $f(q)=m$.
	\end{theorem}
	Let $\overrightarrow{\tilde{q}_s}^{\star}=(1,0,0,0)^\top$ and $\overrightarrow{\tilde{q}_d}^{\star}=\overrightarrow{0}$. Clearly, \((\overrightarrow{\tilde{q}_s}^{\star}, \overrightarrow{\tilde{q}_d}^{\star}) \in \Omega\). Define
	\begin{equation*}
		L_\star := \left\{ (\overrightarrow{\tilde{q}_s},\overrightarrow{\tilde{q}_d}) \;\middle|\; f(\overrightarrow{\tilde{q}_s},\overrightarrow{\tilde{q}_d}) \leq f(\overrightarrow{\tilde{q}_s}^{\star}, \overrightarrow{\tilde{q}_d}^{\star}) \right\}.
	\end{equation*}
	It is easy to verify that \(L_\star\) is nonempty, bounded, and closed. Since \(\Omega\) is closed, we conclude that \(\Omega\cap L_{\star}\) is a closed and bounded subset, and therefore compact.
	As the objective function \(f(\overrightarrow{\tilde{q}_s},\overrightarrow{\tilde{q}_d})\) is continuous, Theorem \ref{EVT} guarantees that the minimum is attained on the compact set \(\Omega\cap L_{\star}\). 
	That is, there exists \((\overrightarrow{\tilde{q}_s}^{\#},\overrightarrow{\tilde{q}_d}^{\#}) \in \Omega\cap L_\star\) such that
	\[
	f(\overrightarrow{\tilde{q}_s}^{\#},\overrightarrow{\tilde{q}_d}^{\#}) = \min_{(\overrightarrow{\tilde{q}_s},\overrightarrow{\tilde{q}_d})\in\Omega} f(\overrightarrow{\tilde{q}_s},\overrightarrow{\tilde{q}_d}),
	\]
	which confirms the existence of a global minimum of problem~\eqref{Objective Function}.

	We define the constraints in problem~\eqref{Objective Function} as follows:
	\[
	c_1(\overrightarrow{\tilde{q}_s}, \overrightarrow{\tilde{q}_d}) := \|\overrightarrow{\tilde{q}_s}\|_2^2 - 1, 
	\quad 
	c_2(\overrightarrow{\tilde{q}_s}, \overrightarrow{\tilde{q}_d}) := \overrightarrow{\tilde{q}_d}^\top \overrightarrow{\tilde{q}_s}.
	\]
	It is evident that the gradients \(\nabla c_1(\overrightarrow{\tilde{q}_s}, \overrightarrow{\tilde{q}_d})\) and \(\nabla c_2(\overrightarrow{\tilde{q}_s}, \overrightarrow{\tilde{q}_d})\) are linearly independent for any \((\overrightarrow{\tilde{q}_s}, \overrightarrow{\tilde{q}_d}) \in \Omega\), which implies that the LICQ holds for problem~\eqref{Objective Function}. By Theorem~\ref{First-Order}, we obtain the following proposition.
	\begin{proposition}\label{proposition1}
		The global minimizer \((\overrightarrow{\tilde{q}_s}^{\#}, \overrightarrow{\tilde{q}_d}^{\#})\) of problem~\eqref{Objective Function} satisfies the KKT conditions.
	\end{proposition}
	Our next objective is to find all feasible pairs \((\overrightarrow{\tilde{q}_s}, \overrightarrow{\tilde{q}_d})\) that satisfy the KKT conditions associated with problem~\eqref{Objective Function}. Among these candidate solutions, the one that attains the smallest value of the objective function corresponds to the desired projection point.
	
	\subsection{The algorithm}\label{Sec 3.3}
	In this subsection, we present our algorithm in a systematic manner. The algorithm is developed based on different cases determined by the relationship between the standard part \(\overrightarrow{\tilde{a}_s}\) and the dual part \(\overrightarrow{\tilde{a}_d}\). Specifically, we consider two scenarios: \(\overrightarrow{\tilde{a}_s} = 0\) and \(\overrightarrow{\tilde{a}_s} \neq 0\). The detailed procedure is provided in Algorithm \ref{alg:projection}.

\begin{figure}[h]
		\centering
		\includegraphics[width=0.9\textwidth]{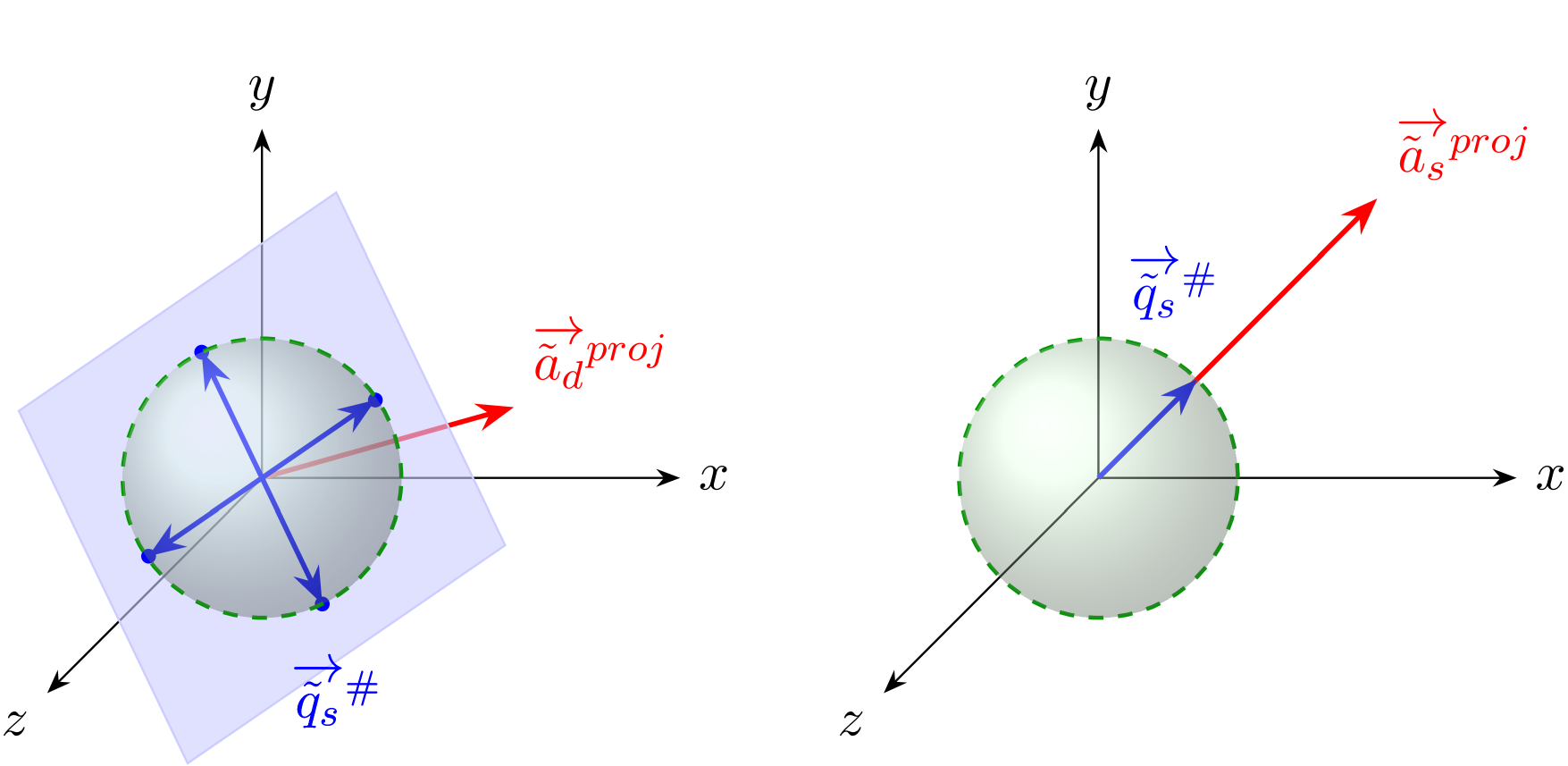}
		\caption{A three-dimensional geometric analogy for the vector projection problem. Left (\textbf{Case 1}): The non-unique solutions $\protect\overrightarrow{\tilde{q}_s}^{\#}$ lie on the blue circular intersection of the green unit sphere and a plane. Right (\textbf{Case 2.1}): The unique solution $\protect\overrightarrow{\tilde{q}_s}^{\#}$ is the radial projection of the input vector $\protect\vec{\tilde{a}}_s$ onto the green unit sphere.}\label{figure1}
	\end{figure}

	\textbf{Case 1: $\overrightarrow{\tilde{a}_s}=\overrightarrow{0}$.} 
	In this case, the constrained optimization problem \eqref{Objective Function} simplifies to 
	\begin{equation*} 
		\min_{\overrightarrow{\tilde{q}_s},\, \overrightarrow{\tilde{q}_d}}  \:\frac{1}{2} \left\| \overrightarrow{\tilde{q}_s}\right\| _2^2+\frac{1}{2} \left\| \overrightarrow{\tilde{q}_d}-\overrightarrow{\tilde{a}_d}\right\| _2^2  \ \ \text{subject to} \ \  \left\| \overrightarrow{\tilde{q}_s}\right\| _2^2 = 1,\; \overrightarrow{\tilde{q}_d}^\top \overrightarrow{\tilde{q}_s} = 0. 
	\end{equation*}
	It can be seen that the globally optimal solutions are given by \(\overrightarrow{\tilde{q}_d}^\# = \overrightarrow{\tilde{a}_d}\), and any \(\overrightarrow{\tilde{q}_s}^\#\) satisfying \(\left\| \overrightarrow{\tilde{q}_s}^\#\right\|  _2^2 = 1\) and $\left(\overrightarrow{\tilde{q}_s}^\# \right)^\top \overrightarrow{\tilde{a}_d}=0 $ is a solution to the problem.  In this case, the optimal solution for $\overrightarrow{\tilde{q}_s}$ is non-unique. It lies on the intersection of the unit hypersphere and the hyperplane orthogonal to the input vector $\overrightarrow{\tilde{a}_d}$ (Fig. \ref{figure1}, left panel).
	%
	
	We now introduce a practical rule for selecting 
	$\overrightarrow{\tilde{q}_s}^\#$ despite this non-uniqueness. 
	Let $\overrightarrow{e_1} = (1, 0, 0, 0)^\top$ and $\overrightarrow{e_2} = (0, 1, 0, 0)^\top$. If \(\overrightarrow{\tilde a_d}=\overrightarrow0\), we set \(\overrightarrow{\tilde q_s}^{\#}=\overrightarrow {e_1}\). If $\overrightarrow{\tilde{a}_{d}}$ is parallel to $\overrightarrow{e_1}$, i.e. $\langle \overrightarrow{\tilde{a}_{d}}, \overrightarrow{e_1} \rangle^2 = \| \overrightarrow{\tilde{a}_{d}} \|_2^2$, we set $\overrightarrow{\tilde{q}_s}^\# = \overrightarrow{e_2}$. Otherwise, we define $\overrightarrow{\tilde{q}_s}^\#$ as the unit projection of $\overrightarrow{e_1}$ onto the orthogonal complement of  $\overrightarrow{\tilde{a}_{d}}$. To summarize, when $\overrightarrow{\tilde{a}_s} = \overrightarrow{0}$, we can compute $\overrightarrow{\tilde{q}_s}^\#$ as follows:
	\begin{equation}\label{deterministic-selection}
		\overrightarrow{\tilde{q}_s}^\# = 	
		\begin{cases}
			\overrightarrow{e_1},& \text{if } \overrightarrow{\tilde a_d}=\overrightarrow 0,\\
			\overrightarrow{e_2},& \text{if } \overrightarrow{\tilde a_d}\neq \overrightarrow 0 \text{ and }\langle \overrightarrow{\tilde{a}_{d}}, \overrightarrow{e_1} \rangle^2 = \| \overrightarrow{\tilde{a}_{d}} \|_2^2,\\
			\frac{\overrightarrow{v}}{\| \overrightarrow{v} \|_2} \ \text{with} \ \ \overrightarrow{v}
			= \overrightarrow{e_1}
			- \frac{\overrightarrow{\tilde{a}_{d}}^\top \overrightarrow{e_1}}
			{\| \overrightarrow{\tilde{a}_{d}} \|_2^2} \, \overrightarrow{\tilde{a}_{d}},
			& \text{otherwise}.
		\end{cases}
	\end{equation}

	\textbf{Case 2: $\overrightarrow{\tilde{a}_s}\neq \overrightarrow{0}$.} In this case, we further distinguish between two subcases: \(\overrightarrow{\tilde{a}_d} = \overrightarrow{0}\) and \(\overrightarrow{\tilde{a}_d} \neq \overrightarrow{0}\).
	
	\textbf{Case 2.1: $\overrightarrow{\tilde{a}_d}=\overrightarrow{0}$.}
	In this case, the constrained optimization problem \eqref{Objective Function} becomes
	\begin{equation}\label{obj-21}
		\min_{\overrightarrow{\tilde{q}_s},\overrightarrow{\tilde{q}_d}} \:\frac{1}{2} \left\| \overrightarrow{\tilde{q}_s} - \overrightarrow{\tilde{a}_s}\right\| _2^2 + \frac{1}{2} \left\| \overrightarrow{\tilde{q}_d}\right\| _2^2  
		\quad \text{subject to} \quad \left\| \overrightarrow{\tilde{q}_s}\right\| _2^2 = 1,\; \overrightarrow{\tilde{q}_d}^\top \overrightarrow{\tilde{q}_s} = 0.
	\end{equation}
	Note that the constraint \(\left\| \overrightarrow{\tilde{q}_s}\right\| _2^2 = 1\) holds, and \(\overrightarrow{\tilde{a}_s}\) is a known constant. Therefore, the objective function in \eqref{obj-21} can be simplified to
	\begin{equation}\label{case2.1-equation}
		\min_{\overrightarrow{\tilde{q}_s},\overrightarrow{\tilde{q}_d}} \: -\left\langle \overrightarrow{\tilde{q}_s}, \overrightarrow{\tilde{a}_s} \right\rangle + \frac{1}{2} \|\overrightarrow{\tilde{q}_d}\|_2^2  
		\quad \text{subject to} \quad \left\| \overrightarrow{\tilde{q}_s}\right\| _2^2 = 1,\; \overrightarrow{\tilde{q}_d}^\top \overrightarrow{\tilde{q}_s} = 0.
	\end{equation}
	
	It is evident that the optimal solution for \(\overrightarrow{\tilde{q}_d}\) is the zero vector, i.e., \(\overrightarrow{\tilde{q}_d}^\# = \overrightarrow{0}\).  
	By the Cauchy-Schwarz inequality, we have
	$
	-\left\langle \overrightarrow{\tilde{q}_s}, \overrightarrow{\tilde{a}_s} \right\rangle \geq -\left\| \overrightarrow{\tilde{q}_s}\right\| _2 \left\| \overrightarrow{\tilde{a}_s}\right\| _2,
	$
	with equality if and only if \(\overrightarrow{\tilde{q}_s}=k\overrightarrow{\tilde{a}_s}\) for some constant $k\geq0$.
	Since \(\left\| \overrightarrow{\tilde{q}_s}\right\| _2 = 1\), it follows that the optimal solution for \(\overrightarrow{\tilde{q}_s}\) is given by
	$
	\overrightarrow{\tilde{q}_s}^\# = \frac{\overrightarrow{\tilde{a}_s}}{\left\| \overrightarrow{\tilde{a}_s}\right\| _2}.
	$
	To summarize, in this case, the unique globally optimal solution is given by  
	\begin{equation*}
		\overrightarrow{q_s}^\# = \frac{\overrightarrow{\tilde{a}_s}}{\left\| \overrightarrow{\tilde{a}_s}\right\| _2} \quad \text{and} \quad \overrightarrow{q_d}^\#= \overrightarrow{0}.
	\end{equation*}
	
	In this case, the unique solution $\overrightarrow{\tilde{q}_s}^\#$ is found by radially projecting the input vector $\overrightarrow{\tilde{a}_s}$ onto the unit hypersphere (Fig. \ref{figure1}, right panel).
	
	\textbf{Case 2.2: $\overrightarrow{\tilde{a}_d}\neq\overrightarrow{0}$.}
	In this case, by Proposition~\ref{proposition1}, we know that the global minimizer \((\overrightarrow{\tilde{q}_s}^{\#}, \overrightarrow{\tilde{q}_d}^{\#})\) of problem~\eqref{Objective Function} satisfies the KKT conditions. Our objective is to identify all feasible pairs \((\overrightarrow{\tilde{q}_s}, \overrightarrow{\tilde{q}_d})\) that satisfy these KKT conditions. The projection point is the one among these candidates that minimizes the objective function.
	
	Using real multipliers $\lambda$ and $\mu$, the Lagrangian function of \eqref{Objective Function} is
	\begin{equation*}
		\mathcal{L}\left( \overrightarrow{\tilde{q}_s},\overrightarrow{\tilde{q}_d},\lambda,\mu\right) =\frac{1}{2} \left\| \overrightarrow{\tilde{q}_s}-\overrightarrow{\tilde{a}_s}\right\| _2^2+\frac{1}{2} \left\| \overrightarrow{\tilde{q}_d}-\overrightarrow{\tilde{a}_d}\right\| _2^2+\lambda\left(\left\| \overrightarrow{\tilde{q}_s}\right\| _2^2-1\right) +\mu\overrightarrow{\tilde{q}_d}^\top \overrightarrow{\tilde{q}_s}.
	\end{equation*}
	According to Theorem~\ref{First-Order}, the first-order necessary optimality conditions for \eqref{Objective Function} are given by
	\begin{align}
		\overrightarrow{\tilde{q}_s} - \overrightarrow{\tilde{a}_s} + 2\lambda \overrightarrow{\tilde{q}_s} + \mu \overrightarrow{\tilde{q}_d} &= \overrightarrow{0}\label{case2.2-8}, \\
		\overrightarrow{\tilde{q}_d} - \overrightarrow{\tilde{a}_d} + \mu \overrightarrow{\tilde{q}_s} &= \overrightarrow{0}\label{case2.2-9}, \\
		\left\| \overrightarrow{\tilde{q}_s}\right\| _2^2 - 1 &= 0\label{case2.2-10}, \\
		\overrightarrow{\tilde{q}_d}^\top \overrightarrow{\tilde{q}_s} &= 0\label{case2.2-11}. 
	\end{align}
	Taking the inner product of \eqref{case2.2-8} with \(\overrightarrow{\tilde{q}_s}\) and using constraints \eqref{case2.2-10} and \eqref{case2.2-11}, we obtain
	\begin{equation}\label{xie-lambda}
		\lambda = \frac{\overrightarrow{\tilde{q}_s}^\top \overrightarrow{\tilde{a}_s} - 1}{2}.
	\end{equation}
	Similarly, taking the inner product of \eqref{case2.2-9} with \(\overrightarrow{\tilde{q}_s}\) yields
	\begin{equation}\label{xie-mu}
		\mu = \overrightarrow{\tilde{q}_s}^\top \overrightarrow{\tilde{a}_d}.
	\end{equation}
	Substituting \(\mu\) back into \eqref{case2.2-9}, we obtain
	\begin{equation*}
		\overrightarrow{\tilde{q}_d} = \overrightarrow{\tilde{a}_d} - \left( \overrightarrow{\tilde{q}_s}^\top \overrightarrow{\tilde{a}_d}\right)  \overrightarrow{\tilde{q}_s}.
	\end{equation*}
	Substituting the expressions for \(\lambda\), \(\mu\), and \(\overrightarrow{\tilde{q}_d}\) into \eqref{case2.2-8} gives
	\begin{equation}\label{case2.2-12}
		\overrightarrow{\tilde{q}_s}^\top \overrightarrow{\tilde{a}_s}\, \overrightarrow{\tilde{q}_s} + \overrightarrow{\tilde{q}_s}^\top \overrightarrow{\tilde{a}_d}\, \overrightarrow{\tilde{a}_d} 
		- \left( \overrightarrow{\tilde{q}_s}^\top \overrightarrow{\tilde{a}_d}\right) ^2 \overrightarrow{\tilde{q}_s} - \overrightarrow{\tilde{a}_s} = \overrightarrow{0}.
	\end{equation}
	Next, we intend to take inner products of \eqref{case2.2-12} with vectors \(\overrightarrow{\tilde{a}_s}\) and \(\overrightarrow{\tilde{a}_d}\), which lead to the following two equations
	\begin{equation}\label{case2.2.13}
		\left( \overrightarrow{\tilde{q}_s}^\top \overrightarrow{\tilde{a}_s}\right) ^2-\left( \overrightarrow{\tilde{q}_s}^\top \overrightarrow{\tilde{a}_d}\right) ^2\left( \overrightarrow{\tilde{q}_s}^\top \overrightarrow{\tilde{a}_s}\right) +\left( \overrightarrow{\tilde{q}_s}^\top \overrightarrow{\tilde{a}_d}\right) \left( \overrightarrow{\tilde{a}_d}^\top \overrightarrow{\tilde{a}_s}\right) -\left\| \overrightarrow{\tilde{a}_s}\right\| _2^2=0
	\end{equation}
	and
	\begin{equation}\label{case2.2.14}
		\left( \overrightarrow{\tilde{q}_s}^\top \overrightarrow{\tilde{a}_s}\right)  \left( \overrightarrow{\tilde{q}_s}^\top \overrightarrow{\tilde{a}_d}\right) +\left( \overrightarrow{\tilde{q}_s}^\top \overrightarrow{\tilde{a}_d}\right) \left\| \overrightarrow{\tilde{a}_d}\right\| _2^2-\left( \overrightarrow{\tilde{q}_s}^\top \overrightarrow{\tilde{a}_d}\right) ^3-\overrightarrow{\tilde{a}_s}^\top \overrightarrow{\tilde{a}_d}=0.
	\end{equation}
	Note that when vectors \(\overrightarrow{\tilde{a}_s}\) and \(\overrightarrow{\tilde{a}_d}\) are linearly dependent, \eqref{case2.2.13} and \eqref{case2.2.14} become equivalent and hence reduce to a single equation, which could simplify the problem. Therefore, we further distinguish between the following two subcases: $\left\langle  \overrightarrow{\tilde{a}_s},\overrightarrow{\tilde{a}_d}\right\rangle ^2=\left\| \overrightarrow{\tilde{a}_s}\right\| _2^2\left\| \overrightarrow{\tilde{a}_d}\right\| _2^2$ and $\left\langle \overrightarrow{\tilde{a}_s},\overrightarrow{\tilde{a}_d}\right\rangle ^2\neq\left\| \overrightarrow{\tilde{a}_s}\right\| _2^2\left\|\overrightarrow{\tilde{a}_d}\right\| _2^2$.
	
	\textbf{Case 2.2.1: $\left\langle \overrightarrow{\tilde{a}_s},\overrightarrow{\tilde{a}_d}\right\rangle ^2=\left\| \overrightarrow{\tilde{a}_s}\right\| _2^2\left\| \overrightarrow{\tilde{a}_d}\right\| _2^2$.}
	Since \(\overrightarrow{\tilde{a}_s} \neq \overrightarrow{0}\) and \(\overrightarrow{\tilde{a}_d} \neq \overrightarrow{0}\), this equality implies that \(\overrightarrow{\tilde{a}_s} = k\overrightarrow{\tilde{a}_d}\) for some scalar \(k \neq 0\).
	By \(\left\|\overrightarrow{\tilde q_s}\right\|_2=1\) and \eqref{xie-mu}, the Cauchy--Schwarz inequality gives $-\left\|\overrightarrow{\tilde a_d}\right\|_2 \leq \mu \leq \left\|\overrightarrow{\tilde a_d}\right\|_2 $.
	Moreover, for any fixed \(\overrightarrow{\tilde q_s}\) satisfying \(\left\|\overrightarrow{\tilde q_s}\right\|_2=1\), the optimal dual part is obtained from \eqref{case2.2-9} as 
	\[ 
	\overrightarrow{\tilde q_d} = \overrightarrow{\tilde a_d} -\mu \overrightarrow{\tilde q_s}. 
	\] 
	Substituting this expression and $\overrightarrow{\tilde{a}_s}=k\overrightarrow{\tilde{a}_d}$ into~\eqref{Objective Function}, we obtain 
	\[ 
	f\left(\overrightarrow{\tilde q_s}, \overrightarrow{\tilde q_d}\right) = \frac12\mu^2 -k\mu+ \frac12 + \frac{k^2}{2} \left\| \overrightarrow{\tilde a_d} \right\|_2^2. 
	\] 
	Therefore, the original optimization problem reduces to 
	\[ 
	\min_{ -\left\|\overrightarrow{\tilde a_d}\right\|_2 \leq \mu \leq \left\|\overrightarrow{\tilde a_d}\right\|_2 } \left\{ \frac12\mu^2-k\mu \right\}. 
	\] 
	Hence, \[ \mu^* = \begin{cases} -\left\|\overrightarrow{\tilde a_d}\right\|_2, & \text{if } k<-\left\|\overrightarrow{\tilde a_d}\right\|_2,\\[2mm] k, & \text{if } -\left\|\overrightarrow{\tilde a_d}\right\|_2 \leq k\leq \left\|\overrightarrow{\tilde a_d}\right\|_2,\\[2mm] \left\|\overrightarrow{\tilde a_d}\right\|_2, & \text{if } k> \left\|\overrightarrow{\tilde a_d}\right\|_2. \end{cases} \] 
	
	After computing \(\mu^*\), we construct
	\(\overrightarrow{\tilde q_s}^{\#}\) by using the orthogonal decomposition with respect to
	\(\overrightarrow{\tilde a_d}\). 
	Choose a unit vector \(\overrightarrow u\) satisfying
	\[
	\overrightarrow u^{\top}\overrightarrow{\tilde a_d}=0,
	\qquad
	\left\|\overrightarrow u\right\|_2=1.
	\]
	Analogous to the rule in \eqref{deterministic-selection} and noting that \(\overrightarrow{\tilde a_d}\neq 0\), we can determine \(\overrightarrow u\) as
	\begin{equation}\label{deterministic-selectionu}
		\overrightarrow u= 	
		\begin{cases} 
			\overrightarrow{e_2},& \langle \overrightarrow{\tilde{a}_{d}}, \overrightarrow{e_1} \rangle^2 = \| \overrightarrow{\tilde{a}_{d}} \|_2^2,\\
			\frac{\overrightarrow{v}}{\| \overrightarrow{v} \|_2} \ \text{with} \ \ \overrightarrow{v}
			= \overrightarrow{e_1}
			- \frac{\overrightarrow{\tilde{a}_{d}}^\top \overrightarrow{e_1}}
			{\| \overrightarrow{\tilde{a}_{d}} \|_2^2} \, \overrightarrow{\tilde{a}_{d}},
			& \text{otherwise}.
		\end{cases}
	\end{equation}
	Since \(\left(\overrightarrow{\tilde q_s}^{\#}\right)^{\top}
	\overrightarrow{\tilde a_d}=\mu^*\) and
	\(\left\|\overrightarrow{\tilde q_s}^{\#}\right\|_2=1\), the component of
	\(\overrightarrow{\tilde q_s}^{\#}\) along
	\(\overrightarrow{\tilde a_d}\) is $\frac{\mu^*}{\left\|\overrightarrow{\tilde a_d}\right\|_2^2}\overrightarrow{\tilde a_d}$,
	and the remaining orthogonal component has norm $\sqrt{1-\frac{(\mu^*)^2}{\left\|\overrightarrow{\tilde a_d}\right\|_2^2}}$.
	Therefore, we set
	\begin{equation*}
		\overrightarrow{\tilde q_s}^{\#} = \frac{\mu^*}{ \left\|\overrightarrow{\tilde a_d}\right\|_2^2 } \overrightarrow{\tilde a_d} + \sqrt{ 1- \frac{(\mu^*)^2}{ \left\|\overrightarrow{\tilde a_d}\right\|_2^2 } } \,\overrightarrow u  \quad \text{and} \quad \overrightarrow{\tilde q_d}^{\#} = \overrightarrow{\tilde a_d} - \mu^* \overrightarrow{\tilde q_s}^{\#}. 
	\end{equation*}
	In particular, in the cases $ \mu^*= \left\|\overrightarrow{\tilde a_d}\right\|_2$  or $\mu^*= -\left\|\overrightarrow{\tilde a_d}\right\|_2$, the square-root term vanishes and the formula reduces to 
	\[ 
	\overrightarrow{\tilde q_s}^{\#} = \frac{\mu^*}{ \left\|\overrightarrow{\tilde a_d}\right\|_2^2 } \overrightarrow{\tilde a_d}, \qquad \overrightarrow{\tilde q_d}^{\#} = \overrightarrow{0}. 
	\]

	Hence, for the case where \(\left\langle  \overrightarrow{\tilde{a}_s}, \overrightarrow{\tilde{a}_d} \right\rangle^2 = \left\| \overrightarrow{\tilde{a}_s}\right\| _2^2 \left\| \overrightarrow{\tilde{a}_d}\right\| _2^2\), the complete procedure for determining the optimal solution \(\left( \overrightarrow{\tilde{q}_s}^\#, \overrightarrow{\tilde{q}_d}^\#\right) \) is formally summarized in Stage~I.
	
	\begin{center}
		{
			\setlength{\fboxrule}{1pt}  
			\setlength{\fboxsep}{10pt} 
			
			\fbox{
				\begin{minipage}{0.9\linewidth} 
					\setstretch{0.4} 
					
					\begin{center}
						\textbf{Stage $\mathrm{I}$}
					\end{center}
					
					\par\vspace{0.5\baselineskip} 
					1: Compute	$k := \frac{\left\langle \overrightarrow{\tilde a_s},\,\overrightarrow{\tilde a_d} \right\rangle}
					{\left\|\overrightarrow{\tilde a_d}\right\|_2^2},$
					which satisfies $\overrightarrow{\tilde a_s}=k\,\overrightarrow{\tilde a_d}$.\\
					\noindent
					
					\noindent
					2: Compute \[ \mu^* = \begin{cases} -\left\|\overrightarrow{\tilde a_d}\right\|_2, & \text{if } k<-\left\|\overrightarrow{\tilde a_d}\right\|_2,\\[1mm] k, & \text{if } -\left\|\overrightarrow{\tilde a_d}\right\|_2 \leq k \leq \left\|\overrightarrow{\tilde a_d}\right\|_2,\\[1mm] \left\|\overrightarrow{\tilde a_d}\right\|_2, & \text{if } k> \left\|\overrightarrow{\tilde a_d}\right\|_2. \end{cases} \] 
					\noindent
					3: Select the vector \(\overrightarrow u\) via \eqref{deterministic-selectionu}.
					
					\noindent
					4: Set \[ \overrightarrow{\tilde q_s}^{\#} = \frac{\mu^*}{ \left\|\overrightarrow{\tilde a_d}\right\|_2^2 } \overrightarrow{\tilde a_d} + \sqrt{ 1- \frac{(\mu^*)^2}{ \left\|\overrightarrow{\tilde a_d}\right\|_2^2 } } \,\overrightarrow u, \] and \[ \overrightarrow{\tilde q_d}^{\#} = \overrightarrow{\tilde a_d} - \mu^* \overrightarrow{\tilde q_s}^{\#}. \]
					
				\end{minipage}
			}
		}
	\end{center}

	\textbf{Case 2.2.2: \(\left\langle  \overrightarrow{\tilde{a}_s}, \overrightarrow{\tilde{a}_d} \right\rangle ^2 \neq \left\| \overrightarrow{\tilde{a}_s}\right\| _2^2 \left\| \overrightarrow{\tilde{a}_d}\right\| _2^2\).}  
	In this case, the strict inequality implies that \(\overrightarrow{\tilde{a}_s}\) and \(\overrightarrow{\tilde{a}_d}\) are linearly independent. A related analysis for this case has also been discussed in \cite{Can2023APL}.
	From \eqref{xie-lambda} and \eqref{xie-mu}, we have
	\[
	\overrightarrow{\tilde{q}_s}^\top \overrightarrow{\tilde{a}_s} = 2\lambda + 1 \ \text{and} \ \overrightarrow{\tilde{q}_s}^\top \overrightarrow{\tilde{a}_d} = \mu.
	\]
	Substituting these expressions into \eqref{case2.2.13} and \eqref{case2.2.14}, we obtain the following coupled equations in \((\lambda, \mu)\)
	\begin{equation}\label{case2.2-13}
		(2\lambda + 1)^2 - (2\lambda + 1)\mu^2 + \mu\,\overrightarrow{\tilde{a}_s}^\top \overrightarrow{\tilde{a}_d} - \left\| \overrightarrow{\tilde{a}_s}\right\| _2^2 = 0,
	\end{equation}
	\begin{equation}\label{case2.2-14}
		(2\lambda + 1)\mu + \mu\left\| \overrightarrow{\tilde{a}_d}\right\| _2^2 - \mu^3 - \overrightarrow{\tilde{a}_s}^\top \overrightarrow{\tilde{a}_d} = 0.
	\end{equation}
	From \eqref{case2.2-14}, we can isolate the term \((2\lambda + 1)\mu\) as
	\[
	(2\lambda + 1)\mu = -\mu \left\| \overrightarrow{\tilde{a}_d}\right\| _2^2 + \mu^3 + \overrightarrow{\tilde{a}_s}^\top \overrightarrow{\tilde{a}_d}.
	\]
	Multiplying equation~\eqref{case2.2-13} by \(\mu^2\) and substituting the above identity into it, we eliminate \(\lambda\) and obtain the following quartic polynomial,
	\begin{equation}\label{case2.2-mu}
		-\left\| \overrightarrow{\tilde{a}_d}\right\| _2^2\,\mu^4 + 2\,\overrightarrow{\tilde{a}_s}^\top \overrightarrow{\tilde{a}_d}\,\mu^3 + \left( \left\| \overrightarrow{\tilde{a}_d}\right\| _2^4 - \left\| \overrightarrow{\tilde{a}_s}\right\| _2^2 \right) \mu^2 - 2\,\overrightarrow{\tilde{a}_s}^\top \overrightarrow{\tilde{a}_d}\, \left\| \overrightarrow{\tilde{a}_d}\right\| _2^2\,\mu + (\overrightarrow{\tilde{a}_s}^\top \overrightarrow{\tilde{a}_d})^2 = 0.
	\end{equation}
	From our derivation we know that if \((\lambda^*,\mu^*)\) satisfies the KKT conditions \eqref{case2.2-8}-\eqref{case2.2-11}, then \((\lambda^*,\mu^*)\) also satisfies equations \eqref{case2.2-13} and \eqref{case2.2-14}. Furthermore, \(\mu^*\) must be a root of the quartic equation \eqref{case2.2-mu}. Consequently, we proceed by finding the roots of \eqref{case2.2-mu}.
	We note that \eqref{case2.2-mu} admits at least one real root, as shown by a sign argument given below. Let \(p(\mu)\) denote the left-hand side of \eqref{case2.2-mu}, i.e.,
	\[
	p(\mu):= -\left\| \overrightarrow{\tilde{a}_d}\right\| _2^2\,\mu^4 + 2\,\overrightarrow{\tilde{a}_s}^\top \overrightarrow{\tilde{a}_d}\,\mu^3 + \left( \left\| \overrightarrow{\tilde{a}_d}\right\| _2^4 - \left\| \overrightarrow{\tilde{a}_s}\right\| _2^2 \right) \mu^2 - 2\,\overrightarrow{\tilde{a}_s}^\top \overrightarrow{\tilde{a}_d}\, \left\| \overrightarrow{\tilde{a}_d}\right\| _2^2\,\mu + (\overrightarrow{\tilde{a}_s}^\top \overrightarrow{\tilde{a}_d})^2.
	\] 
	The leading coefficient of \(p(\mu)\) is \(-\left\| \overrightarrow{\tilde a_d}\right\|_2^2\), which is negative; hence \(p(\mu)\to -\infty\) as \(\mu\to +\infty\). Its constant term equals \(\left( \overrightarrow{\tilde a_s}^{\top}\overrightarrow{\tilde a_d}\right)^2 \ge 0\). If \(\overrightarrow{\tilde a_s}^{\top}\overrightarrow{\tilde a_d}=0\), then \(\mu=0\) is a real root. Otherwise, \(p(0)>0\), and by the intermediate value theorem {\cite[Theorem~5.12]{RudinPMA}}, there exists a real root \(\mu^\star\) such that \(p(\mu^\star)=0\). Therefore, \eqref{case2.2-mu} admits at least one real root.
	All real roots \(\mu\) of the quartic equation~\eqref{case2.2-mu} can be found numerically using a standard polynomial root solver, e.g., in MATLAB. Once a real root \(\mu\) is identified, the corresponding value of \(\lambda\) is obtained using equation~\eqref{case2.2-13} if \(\mu = 0\), or equation~\eqref{case2.2-14} if \(\mu \neq 0\).
	
	Subsequently, from the optimality conditions \eqref{case2.2-8} and \eqref{case2.2-9}, the pair \(\left( \overrightarrow{\tilde{q}_s}, \overrightarrow{\tilde{q}_d}\right) \) can be recovered explicitly as
	\begin{equation}\label{case2.3 qsqd}
		\overrightarrow{\tilde{q}_s} = \frac{\overrightarrow{\tilde{a}_s} - \mu \overrightarrow{\tilde{a}_d}}{2\lambda - \mu^2 + 1}\,\text{ and }\, 
		\overrightarrow{\tilde{q}_d} = \frac{-\mu \overrightarrow{\tilde{a}_s} + (2\lambda + 1) \overrightarrow{\tilde{a}_d}}{2\lambda - \mu^2 + 1}.
	\end{equation}
	We note that  \eqref{case2.3 qsqd} is indeed well-defined, that is  $2\lambda - \mu^2 + 1\neq0$.
	We prove this by contradiction. If $2\lambda - \mu^2 + 1=0$, we obtain $\mu^2\overrightarrow{\tilde q_s}+\mu\overrightarrow{\tilde q_d}=\overrightarrow{\tilde a_s}$ from \eqref{case2.2-8}. Multiplying \eqref{case2.2-9} by $\mu$ yields $\mu^2\overrightarrow{\tilde q_s}+\mu\overrightarrow{\tilde q_d}=\mu\,\overrightarrow{\tilde a_d}$.
	Hence $\overrightarrow{\tilde a_s}=\mu\,\overrightarrow{\tilde a_d}$, which implies
	\(\left\langle  \overrightarrow{\tilde{a}_s}, \overrightarrow{\tilde{a}_d} \right\rangle ^2 = \left\| \overrightarrow{\tilde{a}_s}\right\| _2^2 \left\| \overrightarrow{\tilde{a}_d}\right\| _2^2\). This contradicts the assumption in Case~2.2.2.
	Then, for each real root \(\mu\), we compute the associated \(\left( \lambda, \mu\right) \) and then obtain candidate solutions \(\left( \overrightarrow{\tilde{q}_s}, \overrightarrow{\tilde{q}_d}\right) \) via \eqref{case2.3 qsqd}.
	Among all such candidates, we select the one that yields the minimum value of the original objective function in~\eqref{Objective Function}. 
	In particular,
	in the  case $\mu=0$, which is possible only if $\overrightarrow{\tilde a_s}^{\top}\overrightarrow{\tilde a_d}=0$,  equation~\eqref{case2.2-13} reduces to \((2\lambda + 1)^2 = \left\| \overrightarrow{\tilde{a}_s}\right\| _2^2\), and hence the corresponding values of \(\lambda\) are  
	$
	\lambda =\frac{ -1 \pm \left\| \overrightarrow{\tilde{a}_s}\right\| _2}{2}.
	$
	The two candidates lead to objective function values $\dfrac{1}{2}\left(1 - \left\| \overrightarrow{\tilde{a}_s}\right\| _2\right) ^2$ and $\dfrac{1}{2}\left(1 + \left\|\overrightarrow{\tilde{a}_s}\right\|_2\right)^2$, respectively. Since \(\left(1 - \left\| \overrightarrow{\tilde{a}_s}\right\|_2\right)^2 < \left(1 + \left\|\overrightarrow{\tilde{a}_s}\right\|_2\right)^2\) for all \(\left\|\overrightarrow{\tilde{a}_s}\right\|_2 > 0\), the value  
	$
	\lambda = \frac{ -1 + \left\| \overrightarrow{\tilde{a}_s}\right\| _2}{2}
	$
	yields the smaller objective value and is therefore selected.

	Hence, for the case where \(\left\langle \overrightarrow{\tilde{a}_s}, \overrightarrow{\tilde{a}_d} \right\rangle^2 \neq \left\|\overrightarrow{\tilde{a}_s}\right\|_2^2 \left\|\overrightarrow{\tilde{a}_d}\right\|_2^2\), the complete procedure for determining the optimal solution \(\left(\overrightarrow{\tilde{q}_s}^\#, \overrightarrow{\tilde{q}_d}^\#\right)\) is formally summarized in Stage II. To account for numerical errors in the root computation, let $\eta>0$ denote a prescribed feasibility tolerance.

	\begin{center}
		{ 
			\setlength{\fboxrule}{1pt}  
			\setlength{\fboxsep}{10pt} 
			
			\fbox{
				\begin{minipage}{0.9\linewidth} 
					\setstretch{0.3}
					
					\begin{center}
						\textbf{Stage \(\mathrm{II}\)}
					\end{center}

					\vspace{-\baselineskip}

					\noindent 
					1: Solve \eqref{case2.2-mu} and extract all real roots \(\mu^{(1)}, \ldots, \mu^{(\ell)}\), where \(1 \leq \ell \leq 4\). Set $\mathcal{I}_{\mathrm{feas}}:=\emptyset$.\\
					
					\noindent
					2: \textbf{for} \(i = 1, \ldots, \ell\) \textbf{do}\\
					
					\noindent
					3: \hspace{1.5em} \textbf{if} \(\mu^{(i)} = 0\) \textbf{then}\\
					
					\noindent
					4: \hspace{3.0em} Set \(\overrightarrow{\tilde{q}_s}^{(i)} = \dfrac{\overrightarrow{\tilde{a}_s}}{\left\| \overrightarrow{\tilde{a}_s}\right\|_2}\), \(\overrightarrow{\tilde{q}_d}^{(i)} = \overrightarrow{\tilde{a}_d}\).\\
					
					\noindent
					5: \hspace{1.5em} \textbf{else}\\
					
					\noindent
					6: \hspace{3.0em} Compute \(\lambda^{(i)}\) by solving equation~\eqref{case2.2-14}. Evaluate \(\overrightarrow{\tilde{q}_s}^{\,(i)}\) and \(\overrightarrow{\tilde{q}_d}^{\,(i)}\) using~\eqref{case2.3 qsqd}. 
					
					\noindent
					7: \hspace{1.5em} \textbf{end if}\\

					8: \hspace{1.5em} \textbf{if} $\left|\left\|\overrightarrow{\tilde q_s}^{(i)}\right\|_2^2-1\right|\leq\eta\,\,\text{and}\,\,\left|\left(\overrightarrow{\tilde q_d}^{(i)}\right)^{\top}\overrightarrow{\tilde q_s}^{(i)}\right|\leq\eta$ \textbf{then} \\
					
					9: \hspace{3.0em} Compute
						\[
						f^{(i)}
						=
						\frac{1}{2}
						\left\|
						\overrightarrow{\tilde q_s}^{(i)}
						-
						\overrightarrow{\tilde a_s}
						\right\|_2^2
						+
						\frac{1}{2}
						\left\|
						\overrightarrow{\tilde q_d}^{(i)}
						-
						\overrightarrow{\tilde a_d}
						\right\|_2^2,
						\]
						and set
						$\mathcal{I}_{\mathrm{feas}}
						:=
						\mathcal{I}_{\mathrm{feas}}\cup\{i\}$.\\
					
					10: \hspace{1.5em} \textbf{else}\\
					
					11:  \hspace{3.0em} Discard the candidate
						$\left(
						\overrightarrow{\tilde q_s}^{(i)},
						\overrightarrow{\tilde q_d}^{(i)}
						\right)$.\\
					
					12:  \hspace{1.5em} \textbf{end if}\\
					
					13: \textbf{end for}\\
					
					14: Set
						$\overrightarrow{\tilde q_s}^{\#}
						=
						\overrightarrow{\tilde q_s}^{(i^*)}$
						and
						$\overrightarrow{\tilde q_d}^{\#}
						=
						\overrightarrow{\tilde q_d}^{(i^*)}$,
						where
						$i^*\in\arg\min_{i\in\mathcal{I}_{\mathrm{feas}}}f^{(i)}$.
				\end{minipage}
			}
		}
	\end{center}
	
	The above analysis leads to the following unified procedure for projecting a general dual quaternion onto the unit dual quaternion set, as summarized in Algorithm~\ref{alg:projection}. 
	
	\begin{algorithm}[htpb]
		\caption{Projection onto the unit dual quaternion set.}
		\label{alg:projection}
		\begin{algorithmic}[1]
			\State \textbf{Input:} A general dual quaternion $\hat{a}=\tilde{a}_s+\tilde{a}_d\epsilon$, its vector representations $\overrightarrow{\tilde{a}_s}$ and $\overrightarrow{\tilde{a}_d}$, and a tolerance  $\eta>0$.
			\State \textbf{Output:} Projected vectors \(\overrightarrow{\tilde{q}_s}^\#\), \(\overrightarrow{\tilde{q}_d}^\#\) such that \(\hat{q} = \tilde{q}_s^\# + \tilde{q}_d^\#\epsilon\) is a unit dual quaternion.
			\If {$\overrightarrow{\tilde{a}_s}=\overrightarrow{0}$}
			\State Set $\overrightarrow{\tilde{q}_d}^\#=\overrightarrow{\tilde{a}_d}$. Select $\overrightarrow{\tilde{q}_s}^\#$ such that $\left( \overrightarrow{\tilde{q}_s}^\#\right) ^\top\overrightarrow{\tilde{a}_d}=0$ and $\left\| \overrightarrow{\tilde{q}_s}^\#\right\| _2^2=1$, by
			using the deterministic selection rule \eqref{deterministic-selection};
			\Else
			\If {$\overrightarrow{\tilde{a}_d}=\overrightarrow{0}$}
			\State Set $\overrightarrow{\tilde{q}_s}^\#=\frac{\overrightarrow{\tilde{a}_s}}{\left\| \overrightarrow{\tilde{a}_s}\right\| _2}$ and  $\overrightarrow{\tilde{q}_d}^\#=\overrightarrow{0}$.
			\Else
			\If {$\left\langle \overrightarrow{\tilde{a}_s},\overrightarrow{\tilde{a}_d}\right\rangle ^2=\left\| \overrightarrow{\tilde{a}_s}\right\| _2^2\left\| \overrightarrow{\tilde{a}_d}\right\| _2^2$}
			\State Set $\overrightarrow{\tilde{q}_s}^\#$ and $\overrightarrow{\tilde{q}_d}^\#$ by Stage I.
			\Else
			\State Set $\overrightarrow{\tilde{q}_s}^\#$ and $\overrightarrow{\tilde{q}_d}^\#$ by Stage II.
			\EndIf
			\EndIf
			\EndIf
			\State Projected unit dual quaternion $\hat{q}=\tilde{q}_s^\#+\tilde{q}_d^\#\epsilon$.		
		\end{algorithmic}
	\end{algorithm}

	\begin{remark}\label{remark-2}
		Given $\left( \overrightarrow{\tilde{a}_s},\overrightarrow{\tilde{a}_d}\right) $, existing methods either employ the absolute value definition of dual quaternions \cite{cui2024power} or only consider the configuration corresponding to Case 2.2.2 in our classification \cite{Can2023APL}.
		In contrast, we investigate the projection of the unit dual quaternion set under the $2^R$-norm and provide a systematic analysis covering all possible cases.
		In particular, for Case~2.2.2, we prove that every global minimizer satisfies the KKT system, and we identify the globally optimal solution from among the resulting candidates. Furthermore, when \(\mu^{(i)}=0\), we give explicit formulas for the corresponding candidate solutions.
	\end{remark}

\section{Numerical Examples}\label{SEC4}
In this section, we conduct numerical experiments to evaluate the numerical performance of the unit dual quaternion projection algorithm (Algorithm \ref{alg:projection}) proposed in this paper. All experiments were implemented in {\sc Matlab} for Windows 11 on a LAPTOP PC with an Intel Core Ultra 7 155H @ 1.40 GHz and 32 GB memory. All experimental results were obtained using double-precision floating-point format.

We use the  function {\tt roots()} in {\sc Matlab} to compute the  roots of the quartic equation \eqref{case2.2-mu} in Stage II of Algorithm \ref{alg:projection}.  A root is considered real if \(|\operatorname{Im}(\mu^{(i)})| < 10^{-8}\).
In the numerical implementation, all zero tests and linear-dependence tests in Algorithm~\ref{alg:projection} are performed with a fixed tolerance. The feasibility tolerance in Stage~II is set to $\eta=10^{-8}$ in all numerical experiments.

For Tables~\ref{table2}, \ref{table5}, and~\ref{table6}, each method was warmed up once and then independently executed 20 times on the complete dataset, with the execution order alternated between repetitions to reduce timing bias. For the proposed method, the measured operations included quartic construction, root computation via \texttt{roots()}, candidate recovery, feasibility checking, and objective comparison. Data loading and error computation were excluded for both methods. The reported total runtime is averaged over these repetitions, and the per-sample runtime is obtained by dividing the total runtime by the sample size.
The real-root distribution is reported only for the quartic equation~\eqref{case2.2-mu} used in our method and is independent of the PM comparison.

\subsection{Synthetic data}

First, we verify that Algorithm \ref{alg:projection} can accurately project general dual quaternions onto the unit dual quaternion set with high precision using simulated data. To this end, we apply the following procedure:
\begin{enumerate}
	\item[1)] For a given $n$, the standard part vectors $\overrightarrow{\tilde{a}_{s_i}} \in \mathbb{R}^4$ ($i=1,\dots,n$) are independently generated using {\sc Matlab} syntax {\tt randn(4,n)} and concatenated into matrix $\mathbf{A}_s=\left[ \overrightarrow{\tilde{a}_{s_1}} ,\, \overrightarrow{\tilde{a}_{s_2}},\,\dots,\, \overrightarrow{\tilde{a}_{s_n}}\right]$. Subsequently, we randomly select 10\% of the column vectors $\overrightarrow{\tilde{a}_{s_i}}$ and set them to zero. The purpose of introducing zero entries is to include degenerate configurations	in the synthetic datasets.

	\item[2)] Randomly generate vectors $\mathbf{t}_i\,(i = 1, \dots, n)$ representing translation in 3D space. They are uniformly distributed within the Euclidean cube $[-5,5]^3$ in $\mathbb{R}^3$. Then convert them to pure quaternions in real vector that form the dual parts $\overrightarrow{\tilde{a}_{d_i}}$ of dual quaternions. The $\left\lbrace \overrightarrow{\tilde{a}_{d_i}}\right\rbrace $ are concatenated into matrix $\mathbf{A_d}=\left[ \overrightarrow{\tilde{a}_{d_1}} ,\, \overrightarrow{\tilde{a}_{d_2}},\,\dots,\, \overrightarrow{\tilde{a}_{d_{n}}}\right] \in \mathbb{R}^{4\times n}$. 
	The same random zero-setting procedure is applied to the dual part $\overrightarrow{\tilde{a}_{d_i}}$ to ensure consistency in the construction of the synthetic datasets.
	
	\item[3)] We take $\left( \overrightarrow{\tilde{a}_{s_i}},\overrightarrow{\tilde{a}_{d_i}}\right) $ as inputs, substitute them into Algorithm \ref{alg:projection} to obtain outputs $\left( \overrightarrow{\tilde{q}_{s_i}},\overrightarrow{\tilde{q}_{d_i}}\right) $, and compute:
	\begin{enumerate}
		\item [1)] The norm error, $E_R=\left| \left\| \overrightarrow{\tilde{q}_{s_i}}\right\| _2^2-1\right|$,  measuring the unitarity of the standard part of the unit dual quaternions.
		\item [2)] The orthogonality error, $E_O=\left| \overrightarrow{\tilde{q}_{s_i}} ^\top\overrightarrow{\tilde{q}_{d_i}} \right|$, quantifying the orthogonality between the dual part and the standard part of the unit dual quaternions. 
	\end{enumerate}
	\item[4)] Draw the cumulative distribution function (CDF) plots of norm error and orthogonality error.
	
\end{enumerate}

\begin{figure}[hptb]
	\centering
	\begin{tabular}{cc}
		\includegraphics[width=0.45\linewidth]{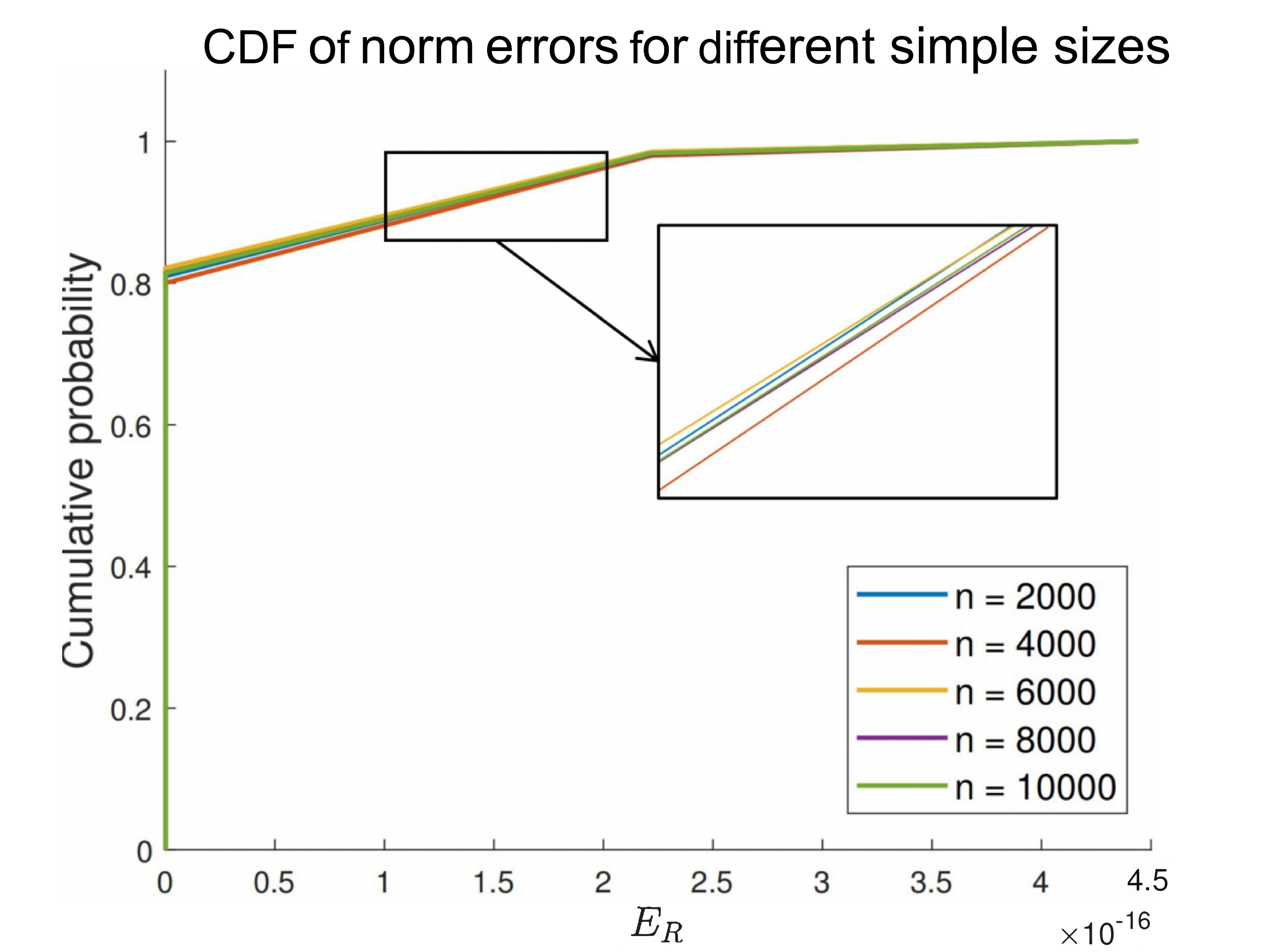}
		\includegraphics[width=0.45\linewidth]{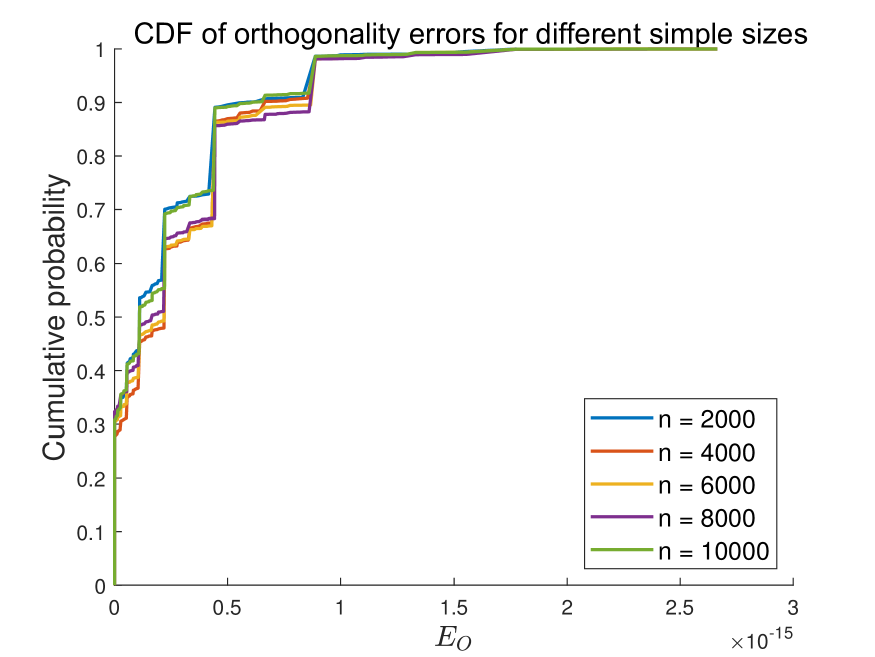}
	\end{tabular}
	\caption{Figures depict the cumulative distribution of norm errors ($E_R$) and orthogonality errors ($E_O$) across different sample sizes}.
	\label{figue2}
\end{figure}

\begin{table}[h]
	\centering
	\caption{Error statistics and OFV comparison for our method and the PM proposed in \cite{cui2024power}.}
	\label{table1}
	\begin{tabular}{cc|ccc} 
		\toprule \tt n & method & $E_R$ & $E_O$ & OFV \\ \midrule \multirow{2}{*}{2000} & Ours & $5.9147\times 10^{-12}$ & $3.2514\times 10^{-13}$ & $8.4081\times 10^{-1}$ \\ & PM & $5.1054\times 10^{-18}$ & $3.1227\times 10^{-15}$ & $2.1440\times 10^{1}$ \\ \midrule \multirow{2}{*}{4000} & Ours & $6.5696\times 10^{-12}$ & $4.0822\times 10^{-13}$ & $8.3853\times 10^{-1}$ \\ & PM & $6.3283\times 10^{-18}$ & $5.0381\times 10^{-15}$ & $4.1974\times 10^{1}$ \\ \midrule \multirow{2}{*}{6000} & Ours & $7.1355\times 10^{-12}$ & $3.8538\times 10^{-13}$ & $8.2032\times 10^{-1}$ \\ & PM & $8.8472\times 10^{-18}$ & $7.5394\times 10^{-15}$ & $7.2294\times 10^{1}$ \\ \midrule \multirow{2}{*}{8000} & Ours & $7.3392\times 10^{-12}$ & $4.0972\times 10^{-13}$ & $8.2081\times 10^{-1}$ \\ & PM & $6.1617\times 10^{-18}$ & $6.2800\times 10^{-15}$ & $5.0087\times 10^{1}$ \\ \midrule \multirow{2}{*}{10000} & Ours & $7.1995\times 10^{-12}$ & $4.0413\times 10^{-13}$ & $8.2349\times 10^{-1}$ \\ & PM & $5.5725\times 10^{-18}$ & $8.6502\times 10^{-15}$ & $6.2595\times 10^{1}$ \\ \bottomrule 
	\end{tabular} 
\end{table}

\begin{table}[h] \centering \caption{Runtime and real-root distribution of \eqref{case2.2-mu} under synthetic data.} \label{table2} 
	\begin{tabular}{cc|cc|ccc} \toprule $n$ & method & \multicolumn{2}{c|}{runtime (s)} & \multicolumn{3}{c}{real-root distribution of \eqref{case2.2-mu} (\%)} \\ \cmidrule(lr){3-4} \cmidrule(lr){5-7} & & per-sample & total & 1 root & 2 roots & $\ge 3$ roots \\ \midrule
		 \multirow{2}{*}{2000} & Ours & $1.17{\times}10^{-5}$ & 0.0233 & 0.00 & 0.62 & 99.38 \\ & PM & $5.19{\times}10^{-7}$ & 0.0010 & -- & -- & -- \\
		\midrule
		\multirow{2}{*}{4000} & Ours & $1.64{\times}10^{-5}$ & 0.0655 & 0.00 & 0.49 & 99.51 \\ & PM & $5.71{\times}10^{-7}$ & 0.0023 & -- & -- & -- \\ \midrule \multirow{2}{*}{6000} & Ours & $1.55{\times}10^{-5}$ & 0.0930 & 0.00 & 0.31 & 99.69 \\ & PM & $5.73{\times}10^{-7}$ & 0.0034 & -- & -- & -- \\ \midrule \multirow{2}{*}{8000} & Ours & $1.76{\times}10^{-5}$ & 0.1404 & 0.00 & 0.37 & 99.63 \\ & PM & $6.37{\times}10^{-7}$ & 0.0051 & -- & -- & -- \\ \midrule \multirow{2}{*}{10000} & Ours & $1.83{\times}10^{-5}$ & 0.1831 & 0.00 & 0.31 & 99.69 \\ & PM & $5.78{\times}10^{-7}$ & 0.0058 & -- & -- & -- \\
		\bottomrule
	\end{tabular}
\end{table}

The CDF representation characterizes the overall distribution of projection errors. It also allows us to compare the constraint residuals across different sample sizes.
In Fig. \ref{figue2}, we plot the CDF of norm errors $E_R$ (left panel) and orthogonality errors $E_O$ (right panel) computed from the projected solutions produced by Algorithm \ref{alg:projection}, where $ n=2000, 4000, \ldots, 10000$. It can be observed from Fig. 2 that both errors remain very small, and the corresponding CDF curves exhibit similar distributions across the tested sample sizes. These results indicate that the proposed method provides consistent numerical satisfaction of the unit-norm and orthogonality constraints as the number of independently processed samples increases.

Table \ref{table1} summarizes the statistical comparison between our method and the projection method (PM) proposed in \cite{cui2024power} on the synthetic datasets. 
Here, we use ``OFV'' to denote the objective function value of problem \eqref{Objective Function}
\[
f\left( \overrightarrow{\tilde{q}_s}, \overrightarrow{\tilde{q}_d}\right) =
\frac{1}{2} \left\| \overrightarrow{\tilde{q}_s} - \overrightarrow{\tilde{a}_s}\right\| _2^2 + \frac{1}{2} \left\| \overrightarrow{\tilde{q}_d} - \overrightarrow{\tilde{a}_d}\right\| _2^2.
\]
By the  definition of the $2^R$-norm, we know that
$
\left\|\hat{a}-\hat{q}\right\|^2_{2^R}
=
\left\|\overrightarrow{\tilde{a}_s}-\overrightarrow{\tilde{q}_s}\right\|_2^2+
\left\|\overrightarrow{\tilde{a}_d}-\overrightarrow{\tilde{q}_d}\right\|_2^2.
$
Hence, it follows that $\left\|\hat{a}-\hat{q}\right\|_{2^R}=\sqrt{2\,f\left( \overrightarrow{\tilde{q}_s}, \overrightarrow{\tilde{q}_d}\right)}$.
Both methods produce very small numerical errors in satisfying the unit-norm and orthogonality constraints. In particular, for the proposed method, the mean norm and orthogonality errors are of orders $10^{-12}$ and $10^{-13}$, respectively. A crucial distinction emerges when evaluating the OFV.

The core advantage of our algorithm is highlighted by the OFV.
Our method consistently produces an optimal projection, with the OFV values ranging from about $0.82$ to $0.84$. In contrast, the PM yields solutions with much larger OFV, ranging from about 21 to 72 in these tests. This significant difference indicates that although the PM can produce points that satisfy the algebraic constraints, it does not generally return the closest point with respect to the OFV considered in this paper. 

Table~\ref{table2} reports the runtime of our method and the PM, together with the real-root distribution of~\eqref{case2.2-mu} for the proposed	method, on the synthetic datasets.
The total runtime of both methods generally increases with the sample size. For the proposed method, the per-sample runtime ranges from $1.17\times10^{-5}$ to $1.83\times10^{-5}$ seconds, whereas that of the PM ranges from $5.19\times10^{-7}$ to $6.37\times10^{-7}$ seconds. Thus, the proposed method operates at the $10^{-5}$-second level per sample, while the PM is consistently faster.
For the quartic equation, no single-real-root cases were observed in the reported synthetic experiments. Two-real-root cases account for between $0.31\%$ and $0.62\%$ of the samples, while cases with at least	three real roots account for between $99.38\%$ and $99.69\%$. All candidates retained by Stage~II satisfied the feasibility tolerances.

\subsection{Real-world data }
In this section, we employ real-world datasets ``{\tt Dataset 1}'' from \cite{sturm2012benchmark} and  ``{\tt Dataset 2}'' from \cite{Furrer2017FSR} for testing Algorithm \ref{alg:projection}. Each dataset comprises multiple sets of poses exhibiting inherent noise characteristics. For the real-world data experiments, we first separate the data into its standard ($\mathbf{A_s}$) and dual ($\mathbf{A_d}$) parts as input to Algorithm \ref{alg:projection}. Tables \ref{table3} and \ref{table4} report the performance of our method and the projection method proposed in \cite{cui2024power} on the ``{\tt Dataset 1}'' and ``{\tt Dataset 2}'' benchmarks, respectively.
\begin{table}[h]
\centering
\caption{Comparison of error statistics and OFV between our method and the PM proposed in \cite{cui2024power} on ``{\tt Dataset 1}''.}
\label{table3}
\setlength{\tabcolsep}{4pt}
\begin{tabular}{c c c|ccc} 
	\toprule Dataset & Scale & Method &$E_R$ & $E_O$ & OFV\\ 
	\midrule \multirow{2}{*}{FB1 XYZ} & \multirow{2}{*}{3000} & Ours & 7.0996 $\times 10^{-17}$ & 7.0470 $\times 10^{-17}$ & 8.4102$\times 10^{-2}$ \\ & & PM & 1.3693$\times 10^{-17}$ & 8.6338$\times 10^{-17}$ & 1.9192$\times 10^{-1}$ \\ 
	\midrule \multirow{2}{*}{FB1 Desk} & \multirow{2}{*}{2335} & Ours& 7.0210 $\times 10^{-17}$ & 4.0210 $\times 10^{-17}$ & 1.4146$\times 10^{-1}$ \\ & & PM & 6.7517$\times 10^{-18}$ & 8.0664$\times 10^{-17}$ & 2.7779$\times 10^{-1}$ \\ 
	\midrule \multirow{2}{*}{FB1 360} & \multirow{2}{*}{{2870}} & Ours & 7.9593 $\times 10^{-17}$ & 2.0816 $\times 10^{-17}$ & 4.9575$\times 10^{-1}$ \\ & & PM & 8.9746$\times 10^{-18}$ & 1.5441$\times 10^{-16}$ & 7.2899$\times 10^{-1}$ \\ 
	\midrule \multirow{2}{*}{FB1 Floor} & \multirow{2}{*}{{4427}} & Ours & 7.9304 $\times 10^{-17}$ & 3.4889 $\times 10^{-17}$ & 4.2758$\times 10^{-1}$ \\ & & PM & 9.1286$\times 10^{-18}$ & 1.0890$\times 10^{-16}$ & 5.2403$\times 10^{-1}$ \\ 
	\midrule \multirow{2}{*}{FB1 Room} & \multirow{2}{*}{{4887}} & Ours & 7.0843 $\times 10^{-17}$ & 3.4023 $\times 10^{-17}$ & 3.4466$\times 10^{-1}$ \\ & & PM & 7.3606$\times 10^{-18}$ & 1.5385$\times 10^{-16}$ & 6.8438$\times 10^{-1}$ \\ 
	\midrule \multirow{2}{*}{FB2 LNL} & \multirow{2}{*}{{6476}} & Ours & 6.9854 $\times 10^{-17}$ & 7.7009 $\times 10^{-17}$ & 4.3129$\times 10^{-1}$ \\ & & PM & 7.7832$\times 10^{-18}$ & 2.5900$\times 10^{-16}$ & 1.1316$\times 10^{0}$ \\ 
	\midrule \multirow{2}{*}{FB2 LWL} & \multirow{2}{*}{{12284}} & Ours & 6.3216 $\times 10^{-17}$ & 7.6228 $\times 10^{-17}$ & 6.8949$\times 10^{-1}$ \\ & & PM & 9.7791$\times 10^{-18}$ & 2.7185$\times 10^{-16}$ & 1.2298$\times 10^{0}$ \\ 
	\midrule \multirow{2}{*}{FB2 XYZ} & \multirow{2}{*}{36820} & Ours & 6.4124 $\times 10^{-17}$ & 2.9762 $\times 10^{-17}$ & 4.5525$\times 10^{-1}$ \\ & & PM & 1.7887$\times 10^{-17}$ & 2.0445$\times 10^{-16}$ & 9.3538$\times 10^{-1}$\\ 
	\midrule \multirow{2}{*}{FB2 Desk} & \multirow{2}{*}{20957} & Ours & 7.2468 $\times 10^{-17}$ & 6.8488 $\times 10^{-17}$ & 3.6474$\times 10^{-1}$ \\ & & PM & 8.3702$\times 10^{-18}$ & 9.0848$\times 10^{-05}$ & 9.8556$\times 10^{-1}$ \\ 
	\midrule \multirow{2}{*}{FB3 NNF} & \multirow{2}{*}{{1580}} & Ours & 5.0030 $\times 10^{-17}$ & 6.4646 $\times 10^{-17}$ & 5.3829$\times 10^{-1}$ \\ & & PM & 0.0000$\times 10^{0\phantom{-}}$ & 2.8261$\times 10^{-16}$ & 1.3698$\times 10^{0}$ \\ 
	\midrule \multirow{2}{*}{FB3 NNN} & \multirow{2}{*}{{3772}} & Ours & 6.1214 $\times 10^{-17}$ & 5.5870 $\times 10^{-17}$ & 7.0256$\times 10^{-1}$ \\ & & PM & 6.8285$\times 10^{-18}$ & 3.2717$\times 10^{-16}$ & 1.5472$\times 10^{0}$ \\ 
	\midrule \multirow{2}{*}{FB3 LO} & \multirow{2}{*}{{8710}} & Ours & 7.1343 $\times 10^{-17}$ & 6.0879 $\times 10^{-17}$ & 4.4680$\times 10^{-1}$ \\ & & PM & 7.6734$\times 10^{-18}$ & 1.3637$\times 10^{-3}$ & 1.0069$\times 10^{0}$ \\ \bottomrule 
\end{tabular} 
\end{table}

\begin{table}[h]
\centering
\caption{Comparison of error statistics and OFV between our method and the PM proposed in \cite{cui2024power} on ``{\tt Dataset 2}''.}
\label{table4}
\setlength{\tabcolsep}{4pt} 
\begin{tabular}{c c c|ccc} 
	\toprule Dataset & Scale & Method & $E_R$ & $E_O$ & OFV \\ \midrule \multirow{2}{*}{FC1 IO} & \multirow{2}{*}{98} & Ours & 7.6252 $\times 10^{-17}$ & 1.4515 $\times 10^{-18}$ & 1.5518$\times 10^{-1}$ \\ & & PM & 9.0630$\times 10^{-18}$ & 3.1464$\times 10^{-17}$ & 1.6195$\times 10^{-1}$ \\ \midrule \multirow{2}{*}{FC1 Vicon} & \multirow{2}{*}{72} & Ours & 8.3923 $\times 10^{-17}$ & 3.5420 $\times 10^{-17}$ & 3.1193$\times 10^{-1}$\\ & & PM & 3.0840$\times 10^{-18}$ & 9.2844$\times 10^{-17}$ & 4.8741$\times 10^{-1}$\\ \midrule \multirow{2}{*}{FC3 IO1} & \multirow{2}{*}{460} & Ours & 7.2824$\times 10^{-17}$ & 9.1283$\times 10^{-18}$ & 7.9963$\times 10^{-2}$ \\ & & PM & 3.8616$\times 10^{-18}$ & 2.4696$\times 10^{-17}$ & 1.2199$\times 10^{-1}$\\ \midrule \multirow{2}{*}{FC3 IO2} & \multirow{2}{*}{949} & Ours & 1.6500$\times 10^{-17}$ & 1.4633$\times 10^{-18}$ & 6.9513$\times 10^{-1}$ \\ & & PM & 4.4456$\times 10^{-18}$ & 1.4601$\times 10^{-16}$ & 7.1098$\times 10^{-1}$ \\ \midrule \multirow{2}{*}{FC3 Vicon} & \multirow{2}{*}{951} & Ours & 3.6697$\times 10^{-17}$ & 3.3153$\times 10^{-17}$ & 4.3042$\times 10^{-1}$ \\ & & PM & 6.5376$\times 10^{-18}$ & 1.4877$\times 10^{-16}$ & 6.8453$\times 10^{-1}$ \\ \midrule \multirow{2}{*}{FC4 IO} & \multirow{2}{*}{1689} & Ours & 5.8651$\times 10^{-17}$ & 1.0850$\times 10^{-17}$ & 1.5836$\times 10^{0}$ \\ & & PM & 3.9440$\times 10^{-18}$ & 1.4176$\times 10^{-4}$ & 1.6883$\times 10^{0}$ \\ \midrule \multirow{2}{*}{FC4 Vicon} & \multirow{2}{*}{1640} & Ours & 1.7929$\times 10^{-17}$ & 2.5880$\times 10^{-17}$ & 7.5452$\times 10^{-1}$ \\ & & PM & 5.2803$\times 10^{-18}$ & 1.2974$\times 10^{-16}$ & 8.6619$\times 10^{-1}$ \\ \midrule \multirow{2}{*}{FC5 IO} & \multirow{2}{*}{1779} & Ours & 5.5682$\times 10^{-17}$ & 2.8513$\times 10^{-17}$ & 1.3876$\times 10^{0}$ \\ & & PM & 5.8663$\times 10^{-18}$ & 3.3610$\times 10^{-17}$ & 1.4531$\times 10^{0}$ \\ \midrule \multirow{2}{*}{FC5 Vicon} & \multirow{2}{*}{1336} & Ours & 3.4348$\times 10^{-17}$ & 3.4293$\times 10^{-17}$ & 8.1666$\times 10^{-1}$ \\ & & PM & 5.1522$\times 10^{-18}$ & 1.2459$\times 10^{-16}$ & 8.4589$\times 10^{-1}$ \\ \bottomrule 
\end{tabular}
\end{table}
Across both sets of experiments, both methods produce very small mean norm errors. Our method also maintains very small orthogonality errors, typically of order $10^{-16}$ or smaller, whereas the PM produces noticeably larger orthogonality errors on several datasets.
Most importantly, the results for the OFV confirm the superior optimality of our approach. Across nearly all datasets, our method achieves a smaller OFV.

\begin{table}[h]
\centering
\caption{Runtime and real-root distribution of \eqref{case2.2-mu} on ``{\tt Dataset 1}''.}
\setlength{\tabcolsep}{4pt} 
\label{table5}
\begin{tabular}{cc|cc|ccc} 
	\toprule
	dataset & method
	& \multicolumn{2}{c}{runtime (s)}
	& \multicolumn{3}{c}{real-root distribution of \eqref{case2.2-mu} (\%)} \\
	\cmidrule(lr){3-4} \cmidrule(lr){5-7}
	& & per-sample & total & 1 root & 2 roots & $\ge 3$ roots \\
	\midrule
	
	\multirow{2}{*}{FB1 XYZ}
	& Ours
	& $1.8338{\times}10^{-5}$
	& $5.5014{\times}10^{-2}$
	& 0.00 & 0.00 & 100.00 \\
	& PM
	& $8.7396{\times}10^{-7}$
	& $2.6219{\times}10^{-3}$
	& -- & -- & -- \\
	\midrule
	
	\multirow{2}{*}{FB1 Desk}
	& Ours
	& $1.6380{\times}10^{-5}$
	& $3.8248{\times}10^{-2}$
	& 0.00
	& 1.03
	& 98.97 \\
	& PM
	& $7.4807{\times}10^{-7}$
	& $1.7467{\times}10^{-3}$
	& -- & -- & -- \\
	\midrule
	
	\multirow{2}{*}{FB1 360}
	& Ours
	& $1.6576{\times}10^{-5}$
	& $4.7572{\times}10^{-2}$
	& 0.00 & 0.00 & 100.00 \\
	& PM
	& $8.3426{\times}10^{-7}$
	& $2.3943{\times}10^{-3}$
	& -- & -- & -- \\
	\midrule
	
	\multirow{2}{*}{FB1 Floor}
	& Ours
	& $1.6254{\times}10^{-5}$
	& $7.1955{\times}10^{-2}$
	& 0.00
	& 43.37
	& 56.63 \\
	& PM
	& $6.8333{\times}10^{-7}$
	& $3.0251{\times}10^{-3}$
	& -- & -- & -- \\
	\midrule
	
	\multirow{2}{*}{FB1 Room}
	& Ours
	& $1.4558{\times}10^{-5}$
	& $7.1146{\times}10^{-2}$
	& 0.00 & 0.00 & 100.00 \\
	& PM
	& $6.7062{\times}10^{-7}$
	& $3.2773{\times}10^{-3}$
	& -- & -- & -- \\
	\midrule
	
	\multirow{2}{*}{FB2 LNL}
	& Ours
	& $1.5145{\times}10^{-5}$
	& $9.8082{\times}10^{-2}$
	& 0.00 & 0.00 & 100.00 \\
	& PM
	& $7.1595{\times}10^{-7}$
	& $4.6365{\times}10^{-3}$
	& -- & -- & -- \\
	\midrule
	
	\multirow{2}{*}{FB2 LWL}
	& Ours
	& $1.8643{\times}10^{-5}$
	& $2.2901{\times}10^{-1}$
	& 0.00 & 0.00 & 100.00 \\
	& PM
	& $8.5113{\times}10^{-7}$
	& $1.0455{\times}10^{-2}$
	& -- & -- & -- \\
	\midrule
	
	\multirow{2}{*}{FB2 XYZ}
	& Ours
	& $1.7106{\times}10^{-5}$
	& $6.2984{\times}10^{-1}$
	& 0.00 & 0.00 & 100.00 \\
	& PM
	& $7.9518{\times}10^{-7}$
	& $2.9279{\times}10^{-2}$
	& -- & -- & -- \\
	\midrule
	
	\multirow{2}{*}{FB2 Desk}
	& Ours
	& $2.4010{\times}10^{-5}$
	& $5.0317{\times}10^{-1}$
	& 0.00 & 0.00 & 100.00 \\
	& PM
	& $9.9667{\times}10^{-7}$
	& $2.0887{\times}10^{-2}$
	& -- & -- & -- \\
	\midrule
	
	\multirow{2}{*}{FB3 NNF}
	& Ours
	& $5.8187{\times}10^{-5}$
	& $9.1935{\times}10^{-2}$
	& 0.00 & 0.00 & 100.00 \\
	& PM
	& $2.0890{\times}10^{-6}$
	& $3.3006{\times}10^{-3}$
	& -- & -- & -- \\
	\midrule
	
	\multirow{2}{*}{FB3 NNN}
	& Ours
	& $5.9821{\times}10^{-5}$
	& $2.2564{\times}10^{-1}$
	& 0.00
	& 2.31
	& 97.69 \\
	& PM
	& $2.1617{\times}10^{-6}$
	& $8.1540{\times}10^{-3}$
	& -- & -- & -- \\
	\midrule
	
	\multirow{2}{*}{FB3 LO}
	& Ours
	& $5.9282{\times}10^{-5}$
	& $5.1635{\times}10^{-1}$
	& 0.00 & 0.00 & 100.00 \\
	& PM
	& $2.2321{\times}10^{-6}$
	& $1.9442{\times}10^{-2}$
	& -- & -- & -- \\
	
	\bottomrule
\end{tabular}
\end{table}
\begin{table}[h]
\centering
\caption{Runtime and real-root distribution of \eqref{case2.2-mu} on ``{\tt Dataset 2}''.}
\setlength{\tabcolsep}{4pt} 
\label{table6}
\begin{tabular}{cc|cc|ccc}
	\toprule
	dataset & Method & \multicolumn{2}{c}{runtime (s)} & \multicolumn{3}{c}{real-root distribution of \eqref{case2.2-mu} (\%)} \\
	\cmidrule(lr){3-4} \cmidrule(lr){5-7}
	& & per-sample & total & 1 root & 2 roots & $\ge 3$ roots \\
	\midrule
	
\multirow{2}{*}{FC1 IO}
& Ours & $2.0104{\times}10^{-5}$ & $1.9702{\times}10^{-3}$ & 0.00 & 100.00 & 0.00 \\
& PM & $1.2514{\times}10^{-6}$ & $1.2264{\times}10^{-4}$ & -- & -- & -- \\
\midrule

\multirow{2}{*}{FC1 Vicon}
& Ours & $1.4737{\times}10^{-5}$ & $1.0610{\times}10^{-3}$ & 0.00 & 100.00 & 0.00 \\
& PM & $8.0000{\times}10^{-7}$ & $5.7600{\times}10^{-5}$ & -- & -- & -- \\
\midrule

\multirow{2}{*}{FC3 IO1}
& Ours & $1.6019{\times}10^{-5}$ & $7.3687{\times}10^{-3}$ & 0.00 & 61.30 & 38.70 \\
& PM & $7.0200{\times}10^{-7}$ & $3.2292{\times}10^{-4}$ & -- & -- & -- \\
\midrule

\multirow{2}{*}{FC3 IO2}
& Ours & $1.5175{\times}10^{-5}$ & $1.4401{\times}10^{-2}$ & 0.00 & 52.37 & 47.63 \\
& PM & $7.0847{\times}10^{-7}$ & $6.7234{\times}10^{-4}$ & -- & -- & -- \\
\midrule

\multirow{2}{*}{FC3 Vicon}
& Ours & $2.0315{\times}10^{-5}$ & $1.9319{\times}10^{-2}$ & 0.00 & 28.60 & 71.40 \\
& PM & $8.5751{\times}10^{-7}$ & $8.1550{\times}10^{-4}$ & -- & -- & -- \\
\midrule

\multirow{2}{*}{FC4 IO}
& Ours & $2.2571{\times}10^{-5}$ & $3.8122{\times}10^{-2}$ & 0.00 & 100.00 & 0.00 \\
& PM & $1.0786{\times}10^{-6}$ & $1.8218{\times}10^{-3}$ & -- & -- & -- \\
\midrule

\multirow{2}{*}{FC4 Vicon}
& Ours & $2.2564{\times}10^{-5}$ & $3.7005{\times}10^{-2}$ & 0.00 & 99.76 & 0.24 \\
& PM & $1.0105{\times}10^{-6}$ & $1.6573{\times}10^{-3}$ & -- & -- & -- \\
\midrule

\multirow{2}{*}{FC5 IO}
& Ours & $2.9684{\times}10^{-5}$ & $5.2808{\times}10^{-2}$ & 0.00 & 100.00 & 0.00 \\
& PM & $1.3679{\times}10^{-6}$ & $2.4334{\times}10^{-3}$ & -- & -- & -- \\
\midrule

\multirow{2}{*}{FC5 Vicon}
& Ours & $4.2565{\times}10^{-5}$ & $5.6866{\times}10^{-2}$ & 0.00 & 92.89 & 7.11 \\
& PM & $1.7228{\times}10^{-6}$ & $2.3016{\times}10^{-3}$ & -- & -- & -- \\

\bottomrule
\end{tabular}
\end{table}
Tables~\ref{table5} and~\ref{table6} report the runtime and real-root distribution of \eqref{case2.2-mu} on  ``{\tt Dataset 1}'' and  ``{\tt Dataset 2}''. For the proposed method, the per-sample runtimes range from $1.4558\times10^{-5}$ to $5.9821\times10^{-5}$ seconds. For the PM, they range from $6.7062\times10^{-7}$ to $2.2321\times10^{-6}$ seconds. Although the PM is consistently faster, the proposed method processes each of the reported real-world datasets in less than one second.
For ``{\tt Dataset 1}'', most samples yield at least three real roots, although two-real-root cases are also observed in FB1 Desk, FB1 Floor, and FB3 NNN. For ``{\tt Dataset 2}'', the proportions of two-real-root and at-least-three-real-root cases vary among the sequences. No single-real-root cases were observed in the reported real-world experiments, and the candidate solutions retained by Stage~II satisfied the prescribed feasibility tolerances.
\section{Conclusion}\label{SEC5}
This paper presents a systematic study on the projection problem for unit dual quaternions, offering a novel algorithm that generalizes existing methods (see Algorithm \ref{alg:projection}). By formulating the problem using the real representation of dual quaternions, we map each dual quaternion to two $4$-dimensional real vectors. Numerical experiments show that the proposed method yields very small constraint residuals and consistently smaller objective values than the PM. Its corrected runtime is at the $10^{-5}$-second level per sample; although slower than the projection method proposed in \cite{cui2024power}, it returns the globally optimal projection for the objective considered.

The results contribute to the research on dual  quaternions, particularly in fields such as robotics, hand-eye calibration, and rigid body motion representation. Future work could explore further optimizations of the algorithm and its applications in more complex robotic systems.

\bibliography{references}
\end{document}